\documentclass[11pt]{article}

\usepackage{amsmath,amssymb,latexsym, color}
\usepackage[normalem]{ulem}
\usepackage{theorem}

\newtheorem{theorem}{Theorem}
\newtheorem{lemma}{Lemma}

\theorembodyfont{\rmfamily} 
\newcommand{\rank}{\operatorname{rank}}

\title{On linearly related orthogonal polynomials in several variables}

\author{Manuel Alfaro$^{a,}$\thanks{Partially supported by MEC of Spain under Grant
MTM2012--36732--C03--02 and Diputaci\'on General de Arag\'on project E--64.}, Ana Pe\~{n}a$^{a,*}$, Teresa E. P\'{e}rez$^{b,}$\thanks{Partially supported by MICINN of Spain
        and by the European Regional Development Fund
        (ERDF) through the grant MTM2011--28952--C02--02, and Junta de
        Andaluc\'{\i}a FQM--0229, P09--FQM--4643 and P11-FQM-7276.}, M. Luisa Rezola$^{a,*}$\\
\\ $^{a}$ {\small Departamento de Matem\'{a}ticas and IUMA,}\\
{Universidad
de Zaragoza (Spain).} \\ $^{b}$ {\small Departamento de
 Matem\'{a}tica Aplicada,}\\{Universidad de Granada (Spain).}}

\date{}
\begin{document}

\maketitle

\bigskip

\begin{abstract}
Let $\{\mathbb{P}_n\}_{n\ge 0}$ and $\{\mathbb{Q}_n\}_{n\ge 0}$ be two monic polynomial systems in several variables satisfying the linear structure relation
$$\mathbb{Q}_n = \mathbb{P}_n + M_n \mathbb{P}_{n-1}, \quad n\ge 1,$$
where $M_n$ are constant matrices of proper size and $\mathbb{Q}_0 = \mathbb{P}_0$. The aim of our work is twofold. First, if both polynomial systems are orthogonal, characterize when that linear structure relation exists in terms of their moment functionals. Second, if one of the two polynomial systems is orthogonal, study when the other one is also orthogonal. Finally, some illustrative examples are presented.

\end{abstract}

\bigskip


\medskip


\newpage

\section{Introduction}

Linear combinations of two families of orthogonal polynomials of one (real or complex) variable have been a subject of great interest for a long time. For instance, it is well known that some families of classical orthogonal polynomials can be expressed as linear combinations of polynomials of the same family with different values of their parameters, the so--called \emph{relations between adjacent families} (e.g. see formulas in Chapter 22 in \cite{AS72} for Jacobi polynomials, or (5.1.13) in \cite{Sz78} for Laguerre polynomials).

The study of such type of linear combinations is related with the concept of \emph{quasi--or\-tho\-go\-na\-li\-ty} introduced by M. Riesz in 1921 (see \cite[p. 64]{Ch78}) as the basis of his analysis of the moment problem.
J. A. Shohat (see \cite{Sh37}) used this notion in connection with some aspects of numerical quadrature; the behaviour of the zeros is also of relevance for problems of approximation theory and interpolation by polynomials, among others.

\bigskip

Likewise, linear combinations of families of multivariate orthogonal polynomials are related with the concept of \emph{quasi--or\-tho\-go\-na\-li\-ty} and they also play an important role in the study of quadrature formulas. Recall the well known results of Gaussian quadrature formulas in the case of one variable (see, e.g. \cite{Ch78}): If  $\{p_n\}$ is a sequence of orthogonal polynomials with respect to either a weight or a definite positive linear functional, then the roots of $p_n+\rho \, p_{n-1}$ with $\rho \in \mathbb R $ are the nodes of a minimal quadrature formula of degree $2n-2$. Moreover, for $\rho=0$ one even obtains a formula of degree $2n-1$. A straightforward extension of these results for higher dimension is not possible. The study of Gaussian cubature started with the classical paper of J. Radon in 1948. The Gaussian cubature formulas of degree $2n-1$ were characterized by Mysovskikh \cite{M81} in terms of the dimension of common zeros of the multivariate orthogonal polynomials. However, these formulas only exist in very special cases and it is the case of degree $2n-2$ that becomes interesting. Here, the linear combinations of multivariate orthogonal polynomials play an important role; again the existence of a Gaussian cubature, now of degree $2n-2$, is given in terms of the dimension of the distinct real common zeros of them, see
 \cite{MP78} and  \cite{S78}. Moreover the nodes of these cubatures formulas are the common zeros of these quasi--orthogonal polynomials. Some progress in this area can be seen in \cite{BSX95}, \cite{SX94}, \cite {X92}, and \cite{X94}.

\bigskip

In recent years there has been a growing interest in linear relations in one variable because of its relationship with several problems, for example:

- The Sobolev orthogonal polynomials, in particular in connection with the notion of coherent pair of measures (\cite{IKNS91}, \cite{MP95}, and \cite{Me97}) and its generalizations.

- The so--called \textit{inverse problem} in the constructive theory of orthogonal polynomials: Given two families of polynomials linearly related, find necessary and sufficient conditions in order to one of them be orthogonal when the other one is orthogonal; see \cite{AMPR03}, \cite{AMPR10}, \cite{AMPR11}, and \cite{MP95}.

- Spectral transformations of moment functionals: Christoffel, Geronimus, Uvarov,...; see \cite{MM92}, \cite{Mar91},  and \cite{Z97}.

- Different properties related to the interlacing of the zeros of particular linear combinations of orthogonal polynomials; see, for instance, \cite{BD05}, \cite{BDR04}, and \cite{DJ09}.

\bigskip

The interest on the orthogonal polynomials of several variables has also increased in recent years. Some problems in which linear relations of multivariate orthogonal polynomials play an important role, are the following: Sobolev orthogonal polynomials (see, e.g. \cite{PPX12} and \cite {X06}),  and the so--called Uvarov and Geronimus modifications of multivariate moment functionals (see, e.g. \cite{DFPP12}, \cite{DFPPX10}, \cite{DGIX09}, and \cite{FPPX10}).

In this context, the \textit{multivariate inverse problem} in the sense described above   appears in a natural way and, as far as we know, it has not been considered in the literature.

Our purpose in this paper is to study  polynomial systems in several variables  $\{\mathbb{P}_n\} _{n \ge 0}$ and $\{\mathbb{Q}_n \}_{n \ge 0}$  satisfying a linear structure relation
\begin{equation*}
\mathbb{Q}_n=\mathbb{P}_n + M_n \mathbb{P}_{n-1}, \quad n\ge 1,
\end{equation*}
where $M_n$ are constant matrices of appropriate size, and $\mathbb{Q}_0=\mathbb{P}_0$. When both polynomial systems are orthogonal, then we prove that only two cases occur, either $M_n\equiv 0, n\ge1$,  or all the matrices $M_n$ have full rank. For these kind of non trivial linear relations, we analyze two inverse problems according to either $\{\mathbb{P}_n\} _{n \ge 0}$ or $\{\mathbb{Q}_n\} _{n \ge 0}$ be orthogonal systems. In the case of one variable the study of these two inverse problems is similar (see \cite{MP95}), however for multivariate orthogonal polynomials, the non--commutativity of the matrices product leads to a quite different situation. Thus, this study is not a simple generalization of the one variable case.

The article is organized as follows. In Section 2 we introduce the basic background that will be needed in the paper. The main results will be stated and developed in Section 3. First part of this section is devoted to study the rank of the matrices $M_n$ in terms of the rank of $M_1$,  when both polynomial systems are orthogonal. Moreover, we give a characterization of the existence  of such linear combination in terms of the relation between the moment functionals. Second part of this section focuses on the study of the multivariate inverse problem. So in Theorems 3 and 4, assuming that one of the polynomial systems is orthogonal we analyze when the other one is also orthogonal. In Section 4 we present a wide set of examples of orthogonal polynomial systems linearly related as above, giving the explicit expressions of the matrices $M_n$. We show particular linear combinations of some bivariate orthogonal polynomial systems introduced by Koornwinder which provide Gaussian cubature formulas of degree $2n-2$ and besides these quasi--orthogonal polynomial systems are also orthogonal.  On the other hand, using the well known Koornwinder's method, we give an example that involves orthogonal polynomials  in two variables on the unit disk. Also we include two examples, namely Krall Laguerre--Laguerre and Krall Jacobi--Jacobi, where the families are orthogonal with respect to quasi--definite moment  functionals. Finally, we deduce {\it relations between adjacent families} of classical orthogonal polynomials in several variables, that is, we express some polynomials as linear combinations of polynomials of the same family with different values of their parameters. In particular we show that these formulas hold for Appell polynomials on the simplex, multiple Jacobi polynomials on the $d$--cube, and multiple Laguerre polynomials on $\mathbb{R}^d_+$. These relations can be seen as a generalization of the ones for Jacobi and Laguerre polynomials in one variable.
\section{Definitions and tools}

Through this paper, we will denote by $\Pi^d$ the
linear space of polynomials in $d$ variables with real coefficients, and
by $\Pi^d_n$ its subspace of polynomials of total degree not greater than
$n$.

Let us denote by ${\mathcal M}_{h\times k}(\mathbb{R})$ the linear space
of $h\times k$ real matrices, and by ${\mathcal
M}_{h\times k}(\Pi^d)$ the linear space of $h\times k$ matrices with polynomial entries. If $h=k$, we will denote
${\mathcal M}_{h\times k} \equiv {\mathcal M}_{h}$, and, in particular, $I_h$
will represent the identity matrix of order $h$. When the dimension of the identity matrix is clear from the context, we will omit the subscript. Given a matrix $M \in
\mathcal{M}_{h}$  we denote by $M^t$ its transpose, and by
$\det(M)$ its determinant. As usual, we say that $M$ is
non--singular if $\det(M)\neq 0$. On the other side, if $M_1, \cdots, M_d$ are matrices of the same size $h\times k$, we define their \emph{joint matrix} $M$ by (\cite[p. 76]{DX01})
\begin{equation*}
M = \begin{pmatrix}
M_1\\
M_2\\
\vdots\\
M_d
\end{pmatrix} = (M_1^t, M_2^t,\cdots, M_d^t)^t,  \qquad M\in {\mathcal M}_{dh\times k}.
\end{equation*}
Next, we will review some basic definitions and properties about multivariate
orthogonal polynomials that we will need along this paper. Most of them can be found in \cite{DX01} which is the main reference in this work.

Let $\mathbb{N}_0$ denote
the set of nonnegative integers. For a multi--index
$\nu=(\nu_1,\dots,\nu_d) \in\mathbb{N}_0^d$, and
$\mathtt{x}=(x_1,\dots,x_d) \in \mathbb{R}^d$ we define a monomial in $d$ variables as
$$\mathtt{x}^{\nu}=
x_1^{\nu_1}\cdots x_d^{\nu_d}.$$ The nonnegative integer $|\nu|=\nu_1+
\dots+\nu_d$ is called the \emph{total degree} of $\mathtt{x}^{\nu}$.

For a fixed total degree $n\ge0$, the cardinal $r^d_n$ of the set of independent monomials
of total degree $n$ is
$$r^d_n = \binom{n+d-1}{d-1}.$$
It is known that there is no natural order for the monomials. In this work, we will use the
\emph{graded lexicographical order}, that is, we order the monomials by the total degree,
and within the monomials of the same total degree, we use the reverse lexicographical order.

For $n\ge 0$, let
$$\{P^n_{\alpha_1}(\mathtt{x}), P^n_{\alpha_2}(\mathtt{x}), \ldots,  P^n_{\alpha_{r_n^d}}(\mathtt{x})\},$$
be $r^d_n$ polynomials of total degree $n$ independent modulus $\Pi^d_{n-1}$,
where $\alpha_1, \alpha_2, \ldots, \alpha_{r_n^d}$ are the elements in
$\{\alpha \in \mathbb{N}^d_0: |\alpha| = n\}$ arranged according to the
reverse lexicographical order. Then we use the column vector notation
$$
\mathbb{P}_n = \mathbb{P}_n(\mathtt{x}) = \begin{pmatrix}
P^n_{\alpha_1}(\mathtt{x})\\ P^n_{\alpha_2}(\mathtt{x})\\ \vdots\\  P^n_{\alpha_{r_n^d}}(\mathtt{x})
\end{pmatrix} =
 (P^n_{\alpha_1}(\mathtt{x}), P^n_{\alpha_2}(\mathtt{x}), \ldots,  P^n_{\alpha_{r_n^d}}(\mathtt{x}))^t.
$$
The sequence of polynomial column vectors $\{\mathbb{P}_n\}_{n\ge0}$ will be called a \emph{polynomial
system} (PS).

Observe that a PS is a sequence of vectors whose dimension and
total degree are increasing: $\mathbb{P}_0$ is a constant,
$\mathbb{P}_1$ is a column vector of dimension $r^d_1$ of multivariate independent
polynomials of total degree 1, $\mathbb{P}_2$ is a column vector
of dimension $r^d_2$ whose elements are multivariate independent polynomials of total
degree 2, and so on.
\noindent
The simplest case of \emph{polynomial system} is the so--called \emph{canonical polynomial
system}, defined as
$$\{\mathbb{X}_n\}_{n\ge0} = \{(\mathtt{x}^{\alpha_1}, \mathtt{x}^{\alpha_2}, \ldots,  \mathtt{x}^{\alpha_{r_n^d}})^t: |\alpha_i| =n\}_{n\ge0}.
$$
Using the vector notation, for a given polynomial system $\{\mathbb{P}_n\}_{n\ge0}$, the vector polynomial
$\mathbb{P}_n$ can be written as
$$\mathbb{P}_n(\mathtt{x}) = G_{n,n}\, \mathbb{X}_n + G_{n,n-1}\, \mathbb{X}_{n-1} + \cdots + G_{n,0}\, \mathbb{X}_0,$$
where $G_n = G_{n,n}$ is called the \emph{leading coefficient} of $\mathbb{P}_n$, which is a square matrix of size $r_n^d$. Moreover, since $\{\mathbb{P}_m\}_{m=0}^n$ form a basis of $\Pi^d_n$, then $G_n$ is invertible.

We will say that two PS $\{\mathbb{P}_n\}_{n\ge0}$ and
$\{\mathbb{Q}_n\}_{n\ge0}$ have the same leading coefficient if $\mathbb{P}_n$ and $\mathbb{Q}_n$ have the same leading coefficient for $n\ge 0$, that is, if the entries of the vector $\mathbb{P}_n
- {\mathbb{Q}}_n$ are polynomials in $\Pi_{n-1}^d$, for $n\ge 1$.

In addition, a polynomial system is called \emph{monic} if every
polynomial contains only one monic term of higher degree, that is, for $n\ge 0$,
$$P^n_{\alpha_k}(\mathtt{x})= \mathtt{x}^{\alpha_{k}} + R(\mathtt{x}), \quad 0 \le k \le r_n^d,$$
where $|\alpha_k| = n$, and $R(\mathtt{x})\in \Pi^d_{n-1}$. Equivalently, a \emph{monic polynomial system} is a polynomial system
such that its leading coefficient is the identity matrix, i. e., $G_n = I_{r_n^d}$, for $n\ge 0$.

Let $s=(s_\alpha)_{\alpha \in \mathbb{N}_0^d}$ be a multi--sequence of real numbers.
We define a linear functional $u$ on $\Pi^d$ by means of the moments
$$\langle u, \mathtt{x}^{\alpha} \rangle= s_\alpha,$$
and extend it by linearity. The linear functional $u$ will be called a \emph{moment functional}.

Recall some operations acting over a moment functional $u$:

\begin{itemize}

\item the {\it action of $u$ over a polynomial matrix}
$$
\langle u, M \rangle :=\left(\langle u, m_{i,j}(\mathtt{x})\rangle
\right)_{i,j=1}^{h,k}\in {\mathcal M}_{h\times k}(\mathbb{R}),$$
where $M=\left(m_{i,j}(\mathtt{x})\right)_{i,j=1}^{h,k} \in {\mathcal M}_{h\times
k}(\Pi^d).$

\item the {\it left product of a polynomial $p\in \Pi^d$ times $u$}
$$\langle p \, u, q \rangle := \langle u, p\, q \rangle,\quad \forall q\in \Pi^d.$$

\item the {\it left product of a matrix of polynomials $M$ times $u$}
$$ \langle M \,u, q \rangle := \langle u,M^t\, q\rangle, \quad \forall
q\in \Pi^d, \quad \forall M\in {\mathcal M}_{h\times
k}(\Pi^d).$$

\item the {\it left product of a matrix of constants $M$ times $u$ acting over a polynomial matrix}
$$ \langle M \,u, N \rangle := \langle u,M^t\, N\rangle = M^t \,\langle u, N\rangle,$$
$\forall M\in {\mathcal M}_{h\times
k}(\mathbb{R}), \quad \forall N\in {\mathcal M}_{h\times
l}(\Pi^d).
$
\end{itemize}

We say that a polynomial $p \in \Pi_n^d$ is {\it orthogonal} with respect
to $u$ if
$$\langle u, p\, q \rangle = 0, \qquad \forall q \in \Pi_{n-1}^d.
$$
The orthogonality can be expressed in terms of a PS $\{\mathbb{P}_n\}_{n\ge 0}$ as
\begin{equation*}
    \langle u, \mathbb{P}_{n}\,\mathbb{P}_{m}^t\rangle  =
    \left\{\begin{array}{ccc} 0 \in {\mathcal M}_{r_n^d\times r_m^d}, & if & n\neq m, \\
    H_n \in {\mathcal M}_{r_n^d\times r_n^d}, & if & n=m,
\end{array}\right.
\end{equation*}
where $H_n$ is a symmetric and non--singular matrix. We shall call
$\{\mathbb{P}_n\}_{n\geq0}$
an \emph{orthogonal polynomial system} (OPS).

A moment functional $u$ is called \emph{quasi--definite} (\cite{DX01}, p. 79) if there is a basis $B$ of
$\Pi^d$ such that for any polynomials $p, q\in B$,
$$\langle u, p\,q\rangle =0, \quad \textrm{if} \quad p\neq q, \quad \textrm{and}\quad
\langle u, p^2\rangle \neq0.$$

The moment functional $u$ is quasi--definite if and only if there exists an
OPS with respect to $u$. If $u$ is quasi--definite, then there
exists a unique \emph{monic orthogonal polynomial system} (MOPS) with respect to $u$.

Moreover, $u$ is \emph{positive definite} if $\langle u, p^2\rangle >0$, for all $p\neq 0$, $p\in \Pi^d$. If $u$ is positive definite, then it is quasi--definite, and it is possible to construct an \emph{orthonormal polynomial system}, that is, an \emph{orthogonal polynomial} system such that $\langle u, \mathbb{P}_{n}\,\mathbb{P}_{n}^t\rangle  = I_{r_n^d}$.

As in the scalar case, orthogonal polynomials in several variables are characterized by a vector--matrix
three term relation (see Theorem 3.2.7 in \cite{DX01}, p. 79). More precisely,
\begin{theorem}[\cite{DX01}] \label{Favard}
Let $\{\mathbb{P}_n\}_{n\ge0} = \{ P^n_\alpha(\mathtt{x}): |\alpha| = n, n\in \mathbb{N}_0\}, \mathbb{P}_0=1$, be an arbitrary sequence in $\Pi^d$. Then the following statement are equivalent.

\begin{enumerate}
\item[(1)] There exists a linear functional $u$ which defines a quasi--definite moment functional on $\Pi^d$ and which makes $\{\mathbb{P}_n\}_{n\ge0}$ an orthogonal basis in $\Pi^d$.

\item[(2)] For $n\ge 0$, $1\le i\le d$, there exist matrices $A_{n,i}$, $B_{n,i}$ and $C_{n,i}$ of respective sizes
$r_n^d\times r_{n+1}^d, \, r_n^d\times r_{n}^d$ and $ r_n^d \times r_{n-1}^d$, such that

\begin{enumerate}

\item the polynomials $\mathbb{P}_n$ satisfy the three term relation
\begin{equation} \label{RR3TP}
 x_i \mathbb{P}_n = A_{n,i} \mathbb{P}_{n+1} + B_{n,i} \mathbb{P}_{n} +C_{n,i} \mathbb{P}_{n-1}, \quad 1\le i\le d,
\end{equation}
with $\mathbb{P}_{-1}=0$ and $C_{-1,i} =0$,

\item for $n\ge 0$ and $1\le i\le d$, the matrices $A_{n,i}$ and $C_{n+1,i}$ satisfy the rank conditions
\begin{equation}\label{condicion  rango 1}
\rank\, A_{n,i} = \rank\,  C_{n+1,i}=r_n^d,
\end{equation}
and, for the joint matrix $A_n$ of $A_{n,i}$, and the joint matrix $C_{n+1}^t$ of $C_{n+1,i}^t$,
\begin{equation}\label{condicion  rango 2}
\rank\, A_{n} = \rank\,  C_{n+1}^t = r_{n+1}^d.
\end{equation}
\end{enumerate}
\end{enumerate}
\end{theorem}

\bigskip

The version of this theorem for orthonormal polynomial systems $\{\mathbb{P}_n\}_{n\ge0}$ is obtained by changing $C_{n+1,i}$ by $A^t_{n,i}$, \,$ 1 \le i \le d, \, n \ge 0 $.

\bigskip

When the orthogonal polynomial system $\{\mathbb{P}_n\}_{n\ge0}$ is monic, comparing the highest coefficient matrices at both sides of \eqref{RR3TP}, it follows that $A_{n,i} = L_{n,i}$, for $n\ge 0$, and $1\le i\le d$, where $L_{n,i}$ are matrices of size $r^d_n\times r^d_{n+1}$  defined by
$$L_{n,i} \, \mathtt{x}^{n+1} = x_i \, \mathtt{x}^n, \quad 1 \leq i \leq d\,.$$ These matrices verify $L_{n,i}L^t_{n,i}=I_{r^d_n}$, and $\rank L_{n,i} = r^d_n$; moreover, the rank of the joint matrix $L_n$ of $L_{n,i}$ is $r_{n+1}^d$ (\cite[p. 77]{DX01}).

For the particular case $d=2$, we have that $L_{n,i}$, $i=1,2$, are the $(n+1)\times(n+2)$ matrices defined as
$$L_{n,1} = \begin{pmatrix}
1 & 0 & \cdots & 0 & 0\\
0 & 1 &  \cdots & 0 & 0\\
\vdots & \vdots & \ddots & \vdots & \vdots\\
0 & 0 & \cdots & 1 & 0\\
\end{pmatrix}_{(n+1)\times (n+2)} \quad \text{and} \qquad \quad
L_{n,2} = \begin{pmatrix}
0 & 1 &  0 & \cdots & 0 \\
0 & 0 &  1 & \cdots & 0 \\
\vdots & \vdots & \vdots & \ddots & \vdots\\
0 & 0 & 0& \cdots & 1\end{pmatrix}_{(n+1)\times (n+2)}.$$

\section{Main results}

In this section we consider two monic polynomial systems $\{\mathbb{P}_n \}_{n \ge 0}$ and $\{\mathbb{Q}_n \}_{n \ge 0}$ related by
\begin{equation}\label{relacion 1-2}
 \mathbb{Q}_n =\mathbb{P}_n + M_n \mathbb{P}_{n-1}, \quad {n\ge 0,}
\end{equation}
where $M_n\in {\cal{M}}_{r^d_{n}\times r^d_{n-1}}$, $n\ge 1$, are constant matrices and $\mathbb{Q}_0=\mathbb{P}_0$. For convenience, through the paper, we adopt the convention $M_0 \equiv 0$.
From now on, we will say that $\{\mathbb{P}_n \}_{n \ge 0}$ and $\{\mathbb{Q}_n \}_{n \ge 0}$ are \emph{linearly related} by means of \eqref{relacion 1-2}.

The monic character of the polynomial systems in \eqref{relacion 1-2} is superfluous. In fact, for $n\ge 0$, let $E_n$, $F_n$ be non--singular matrices of size $r^d_n$, and define the new polynomial systems $\{\hat{\mathbb{P}}_n\}_{n\ge0}$, $\{\hat{\mathbb{Q}}_n\}_{n\ge0}$ by means of
$$\hat{\mathbb{P}}_n = E_n\, \mathbb{P}_n, \qquad \hat{\mathbb{Q}}_n = F_n\, \mathbb{Q}_n, \quad n\ge 0.
$$
Since $\mathbb{P}_n$ and $\mathbb{Q}_n$ are monic, then $E_n$ and $F_n$ are the leading coefficients of $\hat{\mathbb{P}}_n$ and $\hat{\mathbb{Q}}_n$, respectively.
\newline
Multiplying \eqref{relacion 1-2} by $F_n$, we get
\begin{eqnarray*}
\hat{\mathbb{Q}}_n &=& F_n\, \mathbb{Q}_n = F_n\, \mathbb{P}_n + F_n\,M_n\, \mathbb{P}_{n-1}\\
&=& F_n\,E_n^{-1}\, E_n \mathbb{P}_n + F_n\,M_n\,E^{-1}_{n-1}\,E_{n-1}\, \mathbb{P}_{n-1}\\
&=& F_n\,E_n^{-1}\, \hat{\mathbb{P}}_n + F_n\,M_n\,E^{-1}_{n-1}\,\hat{\mathbb{P}}_{n-1}.
\end{eqnarray*}
Then
\begin{equation}\label{rel1-2non-monic}
\hat{\mathbb{Q}}_n = \hat{K}_n\, \hat{\mathbb{P}}_n + \hat{M}_n\,\hat{\mathbb{P}}_{n-1},
\end{equation}
that is, $\{\hat{\mathbb{P}}_n\}_{n\ge0}$ and $\{\hat{\mathbb{Q}}_n\}_{n\ge0}$ are linearly related by the above expression, where $\hat{K}_n = F_n\,E_n^{-1}$ and $\hat{M}_n = F_n\,M_n\,E^{-1}_{n-1}$. When both polynomial systems have the same leading coefficients then $\hat{K}_n=I_{r_n^d}$, but, in general, $\hat{K}_n$ is a non--singular matrix. Moreover, $\rank \hat{M}_n = \rank M_n$, $n\ge0$, since the rank is unchanged upon left or right multiplication by a non--singular matrix (\cite[p. 13]{HJ85}).

Since the definition of \emph{linearly related} does not depend on particular bases, it is often more convenient to work with monic polynomial systems.

\bigskip

First of all, we analyze the case when both monic polynomial systems are orthogonal, and we deduce some properties about the rank of all matrices $M_n$ in relation (\ref{relacion 1-2}) in terms of the rank of $M_1$.

\begin{lemma}\label{lema rango}
Let $\{\mathbb{P}_n \}_{n \ge 0}$ and $\{\mathbb{Q}_n \}_{n \ge 0}$  be two monic orthogonal polynomial systems linearly related by (\ref{relacion 1-2}). Then

\noindent (i) If $\rank M_1=0$, then $\rank M_n=0$ for every $n \ge 1$,

\noindent (ii) If $\rank M_1=1$, then $\rank M_n=r_{n-1}^d$ for every $n\ge1$.
\end{lemma}

\noindent \textbf{Proof.}
Because of
$$\{\mathbb{Q}_m: m\geq 0\} = \{Q^m_\beta(\mathtt{x}): \,\,\beta =(b_1,\ldots,b_d)\in \mathbb{N}^d_0,\,\, |\beta| =m, \, \,m\ge 0\},
$$
is an algebraic
basis in $\Pi^d$, we can associate to it the corresponding dual basis on the algebraic dual space of $\Pi^d$
$$\{f^n_\alpha: \,\,\alpha=(a_1,\ldots,a_d)\in \mathbb{N}^d_0, \,\,|\alpha| =n, \,\,n\ge 0\},
$$
where $f^n_\alpha$ is the linear functional defined as
$$\langle f^n_\alpha, Q^m_\beta\rangle = \delta_{n,m}\delta_{a_1,b_1}\cdots \delta_{a_d,b_d}.
$$
If $f$ is an arbitrary linear functional in the dual space of $\Pi^d$, then it can be written as a linear combination of the elements of the basis, that is,
\begin{equation}\label{fou}
f=\sum_{n=0}^{+\infty} \sum_{|\alpha| = n} \varepsilon^n_\alpha \, f^n_\alpha, \qquad \textrm{where} \quad \varepsilon^n_\alpha = \langle f, Q^n_\alpha\rangle.
\end{equation}
This dual basis can be written as a sequence of row vectors of functionals
$$\Phi_n = (f^n_{\alpha_1}, \ldots, f^n_{\alpha_{r^d_{n}}})_{1\times r^d_{n}},\quad n\ge 0,$$
where $\alpha_1, \ldots, \alpha_{r_n^d}$ are the elements in
$\{\alpha \in \mathbb{N}^d_0: |\alpha| = n\}$ arranged according to the
reverse lexicographical order. Obviously, we can express the duality as follows
$$\langle \Phi_n, \mathbb{Q}^t_m\rangle = \begin{cases}
                                              0 \in {\mathcal M}_{r_n^d\times r_m^d}, & n\neq m,\cr
                                              \cr
                                               I_{r_n^d}, & n=m,\end{cases}$$
and expression \eqref{fou} can be written in a vector form as
$$f=\sum_{n=0}^{+\infty} \Phi_n \, E_n, \qquad E_n =\langle f, \mathbb{Q}_n \rangle\in {\cal{M}}_{r^d_{n}\times 1}.$$
Let $u$ and $v$ be the respective quasi--definite moment functionals associated with the orthogonal polynomial systems $\{\mathbb{P}_n\}_{n\geq 0}$ and $\{\mathbb{Q}_n\}_{n\geq 0}$. Since $\tilde{H}_n=\langle v,  \mathbb{Q}_n\,\mathbb{Q}_n^t\rangle$ is an invertible matrix for $n\ge 0$, we deduce
\begin{eqnarray*}
\langle \mathbb{Q}_n^t\, \tilde{H}_n^{-1}\,v, \mathbb{Q}_m^t\rangle &=&\langle \tilde{H}_n^{-1}\,v,  \mathbb{Q}_n\,\mathbb{Q}_m^t\rangle \\
&=& \tilde{H}_n^{-1}\,\langle v,  \mathbb{Q}_n\,\mathbb{Q}_m^t\rangle =\begin{cases}
                                              0 \in {\mathcal M}_{r_n^d\times r_m^d}, & n\neq m,\cr
                                              \cr
                                               I_{r_n^d}, & n=m,\end{cases}
\end{eqnarray*}
that is, the row linear functionals $\Phi_n$ and $\mathbb{Q}_n^t\, \tilde{H}_n^{-1}\,v$ coincide over the basis
$\{\mathbb{Q}_n\}_{n\geq 0}$, and then $\Phi_n=\mathbb{Q}_n^t\, \tilde{H}_n^{-1}\,v$. Thus, there exist column vectors of constants $E_n\in {\cal{M}}_{r^d_{n}\times 1}$  for $n\ge 0$, such that we can express   $u$ in terms of this dual basis as
\begin{equation*}
u=\sum_{n=0}^{+\infty} \mathbb{Q}_n^t\, \tilde{H}_n^{-1}\, E_n\,v.
\end{equation*}
Observe that
$$
\langle u, \mathbb{Q}_k^t \rangle  =\sum_{n=0}^{+\infty}  E_n^t \, \tilde{H}_n^{-1} \, \langle  v, \mathbb{Q}_n \, \mathbb{Q}^ t_k \rangle = E_k^t.
$$
Now, taking into account relation (\ref{relacion 1-2}), we have
$$\langle u, \mathbb{Q}_k^t \rangle = \langle u, \mathbb{P}_k^t + \mathbb{P}_{k-1}^t\,
M^t_k\rangle =0,\qquad k\ge 2,
$$
and then
\begin{equation*}
u=\big(\mathbb{Q}_1^t \, \tilde{H}_1^ {-1}\,E_1  + \mathbb{Q}_0^t \, \tilde{H}_0^ {-1}\, E_0\big)\,v,
\end{equation*}
where $E_0^t = \langle u, \mathbb{Q}_0^t\rangle = H_0$, and $E_1^t= \langle u, \mathbb{Q}_1^t\rangle = \langle u, \mathbb{P}_1^t + \mathbb{P}_{0}^t\,
M^t_1\rangle = H_0\, M_1^t$.

Therefore, we can write
\begin{equation}\label{primera relacion funcionales}
u = \big(\mathbb{Q}_1^t \, \tilde{H}_1^{-1}\,M_1 + \tilde{H}_0^{-1} \big)H_0\, v = \lambda(\mathtt{x})\,v.
\end{equation}

\noindent
(i) If $\rank M_1=0$, that is $M_1\equiv 0$, then by (\ref{primera relacion funcionales}), $u = ( \tilde{H}_0^{-1}\,H_0 )\, v$. Thus
$\mathbb{P}_n=\mathbb{Q}_n$,  for all $n\ge0$, and again from (\ref{relacion 1-2}) we obtain $M_n\equiv0$.

\noindent
(ii) If $\rank M_1=1$, that is $M_1$ has full rank, then from (\ref{primera relacion funcionales}), $u=\lambda(\mathtt{x})v$ where $\lambda(\mathtt{x})$ is a polynomial of exact total degree one, namely
$$\lambda(\mathtt{x}) = \sum_{i=1}^d a_i\,x_i + b, \quad \mathrm{with} \quad \sum_{i=1}^d|a_i|\neq 0.$$
Using \eqref{relacion 1-2} and the three term relation \eqref{RR3TP}, we get
\begin{eqnarray*}
M_n\, H_{n-1} &=& M_{n}\,\langle u, \mathbb{P}_{n-1}\mathbb{P}_{n-1}^t\rangle = \langle u,\mathbb{Q}_n \,
\mathbb{P}_{n-1}^t\rangle =  \langle \lambda(\mathtt{x}) v,\mathbb{Q}_n \mathbb{P}_{n-1}^t\rangle \\
&=& \sum_{i=1}^d a_i \langle v,\mathbb{Q}_n \,x_i\,
\mathbb{P}_{n-1}^t\rangle = \sum_{i=1}^d a_i\langle v,\mathbb{Q}_n \,\mathbb{P}_n^t\rangle L_{n-1,i}^t = \tilde{H}_n\,\sum_{i=1}^d a_i \, L_{n-1,i}^t,
\end{eqnarray*}
in summary,
\begin{equation}\label{M-H}
M_{n} \,H_{n-1} = \tilde{H}_n\,\left(\sum_{i=1}^d a_i \, L_{n-1,i}^t\right).
\end{equation}
The special shape of the matrices $L_{n-1,i}$ described in the above section, allows to deduce that the rank of
the matrix $\left(\sum_{i=1}^d a_i \, L_{n-1,i}^t\right)$ is $r_{n-1}^d$. Then
$$\rank \,M_n = \rank \,M_n\, H_{n-1} = \rank \tilde{H}_n\,\left(\sum_{i=1}^d a_i \, L_{n-1,i}^t\right) = \rank \,\Big(\sum_{i=1}^d a_i \, L_{n-1,i}^t\Big) = r_{n-1}^d,$$
since $\tilde{H}_n$ and $H_{n-1}$ are non--singular matrices and the rank condition is invariant through non--singular matrices (\cite{HJ85}, p. 13). \quad $\square$

\bigskip

Next, we characterize when two monic orthogonal polynomial systems  $\{\mathbb{P}_n\}_{n \ge 0}$ and $\{\mathbb{Q}_n \}_{n \ge 0}$ are related by a formula as (\ref{relacion 1-2}) in terms of the relation between  their  respective moment functionals.

\begin{theorem} \label{teorema principal} Let $\{\mathbb{P}_n \}_{n \ge 0}$ and $\{\mathbb{Q}_n \}_{n \ge 0}$ be  two monic orthogonal polynomial systems, and let $u$ and $v$ be their quasi--definite moment functionals, respectively. Then the following conditions are equivalent:

\noindent (i) There exist real matrices $M_n \in {\cal{M}}_{r^d_{n}\times r^d_{n-1}}$ with $M_{1}\not\equiv 0$,
such that $\{\mathbb{P}_n \}_{n\ge0}$ and $\{\mathbb{Q}_n \}_{n\ge0}$ are related by
(\ref{relacion 1-2}).

\noindent (ii) There exists a polynomial $\lambda(\mathtt{x})$ of degree one such that
\begin{equation*}
u=\lambda(\mathtt{x})\, v,
\end{equation*}
and $\mathbb{P}_1 \not=\mathbb{Q}_1$.
\end{theorem}
\noindent \textbf{Proof.}

\noindent
$(i)\Rightarrow (ii)$ From (\ref{relacion 1-2}) and  $M_{1}\not\equiv 0$ we have $\mathbb{P}_1 \not=\mathbb{Q}_1$, and from (\ref{primera relacion funcionales}) there exists a polynomial $\lambda(\mathtt{x})$ of  degree one such that $u=\lambda(\mathtt{x})v$ where
$$\lambda(\mathtt{x})=\big(\mathbb{Q}_1^t \, \tilde{H}_1^{-1}\,M_1 + \tilde{H}_0^{-1} \big)H_0.$$

\noindent
$(ii)\Rightarrow (i)$  Consider the Fourier expansion of $\mathbb{Q}_n$ in terms of the
polynomials $\mathbb{P}_n$,
$$\mathbb{Q}_n =\mathbb{P}_n+\sum_{j=0}^{n-1}M_{n,j}\,\mathbb{P}_j.$$
Then
\begin{equation*}
M_{n,j}  = \langle u,\mathbb{Q}_n \mathbb{P}_j^t\rangle \, H_j^{-1}= \langle \lambda(\mathtt{x})\, v, \mathbb{Q}_n \mathbb{P}_j^t\rangle\, H_j^{-1}
 = \langle v,  \mathbb{Q}_n\,\lambda(\mathtt{x})\, \mathbb{P}_j^t\rangle\, H_j^{-1} = 0, \quad  0\le j \le n-2.
\end {equation*}
 Thus
\begin{equation*}
\mathbb{Q}_n=\mathbb{P}_n + M_n \mathbb{P}_{n-1}, \quad n\ge 1,
\end{equation*}
where we denote $M_n \equiv M_{n,n-1}$.

Observe that if the explicit expression for the polynomial $\lambda$ is $\lambda(\mathtt{x}) = \sum_{i=1}^d a_i\,x_i + b$, with $\sum_{i=1}^d|a_i|\neq 0$, as in Lemma \ref{lema rango}, we get formula \eqref{M-H}.

Thus, for $n=1$, it follows
$$M_1 = \tilde{H}_1 \begin{pmatrix}
a_1\\ a_2\\ \vdots\\ a_d\end{pmatrix} \,H_0^{-1} \not \equiv 0. \quad \square$$

\bigskip

In the sequel, we will consider two polynomial systems $\{\mathbb{P}_n \}_{n \ge 0}$ and $\{\mathbb{Q}_n \}_{n \ge 0}$, and we will use the three term relation \eqref{RR3TP} in the adequate conditions. Whenever the system $\{\mathbb{Q}_n \}_{n \ge 0}$ be orthogonal in its corresponding three term relation, we will use the \emph{tilde} notation:
\begin{equation}\label{RR3TQ}
x_i \,\mathbb{Q}_n = \tilde{A}_{n,i} \,\mathbb{Q}_{n+1} + \tilde{B}_{n,i}\, \mathbb{Q}_{n} +\tilde{C}_{n,i}\, \mathbb{Q}_{n-1}, \quad n\ge 0,
\end{equation}
and assume that conditions \eqref{condicion  rango 1} and \eqref{condicion  rango 2}
for $\tilde{A}_{n,i}$ and $\tilde{C}_{n,i}$ are satisfied.

\bigskip

Now, let us analyze the following problem: assuming one of the two polynomial systems related by (\ref{relacion 1-2}) is orthogonal, characterize when the other is also orthogonal. As a consequence of Lemma \ref{lema rango} the case when the matrix $M_1$ in the relation (\ref{relacion 1-2}) has full rank is the only case to be considered, so the condition $\rank M_1 = 1$ will be imposed in what follows.

We obtain the following two characterizations.

\begin{theorem}\label{caracterizacion ortogonalidad 2} Let $\{\mathbb{Q}_n \}_{n \ge 0}$ be
a system of monic orthogonal polynomials satisfying \eqref{RR3TQ}. {Define} recursively a polynomial system $\{\mathbb{P}_n \}_{n \ge 0}$ by (\ref{relacion 1-2}) with $\rank M_1=1$. Then $\{\mathbb{P}_n \}_{n \ge 0}$ is a {monic orthogonal polynomial system} satisfying the three term relation \eqref{RR3TP} if and only if
\begin{equation}\label{formula caracterizacion2}
M_n \, C_{n-1,i}=\tilde{C}_{n,i} \, M_{n-1}, \quad n\ge 2,
\end{equation}
and
\begin{eqnarray}
& ~ & {A}_{n,i}=\tilde{A}_{n,i} = L_{n,i},\label{A}\\
& ~ & {B}_{n,i}=\tilde{B}_{n,i} - M_n A_{n-1,i}+ \tilde{A}_{n,i} M_{n+1},  \label{B}\\
& ~ & {C}_{n,i}=\tilde{C}_{n,i} - M_n B_{n-1,i}+ \tilde{B}_{n,i} M_{n}. \label{C}
\end{eqnarray}
\end{theorem}

\noindent
\textbf{Proof.} Inserting (\ref{relacion 1-2}) in (\ref{RR3TQ}), we have for every $i\in \{1, \dots, d\}$, and $n\ge 1$,
\begin{align*}
&\big(x_i \, I -\tilde{B}_{n,i}\big)\big[\mathbb{P}_n + M_n \mathbb{P}_{n-1}\big]\\
\nonumber & -\tilde{A}_{n,i}
\big[\mathbb{P}_{n+1} + M_{n+1} \mathbb{P}_{n}\big]- \tilde{C}_{n,i}\big[\mathbb{P}_{n-1} + M_{n-1} \mathbb{P}_{n-2}\big]=0.
\end{align*}

Assume that $\{\mathbb{P}_n\}_{n\ge0}$ is an OPS. Then by (\ref{RR3TP}), we get
\begin{align*}
&\big( A_{n,i}- \tilde{A}_{n,i}\big)\mathbb{P}_{n+1}
+ \big( {B}_{n,i} -   \tilde{B}_{n,i} + M_n A_{n-1,i}- \tilde{A}_{n,i} M_{n+1} \big)\mathbb{P}_{n} \\
&+\big( C_{n,i} - \tilde{C}_{n,i} + M_n B_{n-1,i}- \tilde{B}_{n,i} M_{n}   \big)\mathbb{P}_{n-1}\\
&+\big( M_n \, C_{n-1,i} -  \tilde{C}_{n,i} \, M_{n-1}  \big)\mathbb{P}_{n-2} =0.
\end{align*}
Using the fact that $\{\mathbb{P}_n \}_{n \ge 0}$ is a basis of $\Pi^d$ we obtain (\ref{formula caracterizacion2})--(\ref{C}).

\bigskip

Conversely, first of all, we are going to use an induction procedure to verify that $\{\mathbb{P}_n\}_{n\ge0}$ satisfies a three term relation as \eqref{RR3TP}. Take the matrices ${A}_{n,i}, \, {B}_{n,i}$,  and ${C}_{n,i}$ given by (\ref{A}), (\ref{B}), and (\ref{C}), respectively. Multiplying \eqref{relacion 1-2} for $n=1$ by $\tilde{A}_{0,i}$, it is easy to see that
\begin{equation*}
x_i \,\mathbb{P}_0 =\tilde{A}_{0,i}\, \mathbb{P}_{1} + \big(\tilde{B}_{0,i}+ \tilde{A}_{0,i} M_1\big) \mathbb{P}_{0},
\end{equation*}
and so  the first step for the induction procedure is obtained.
Now, we suppose that (\ref{RR3TP}) holds for  $n-1$ and we are going to prove it for $n$.

Multiplying $\mathbb{P}_{n+1}$  by ${A}_{n,i}$ in the relation given in (\ref{relacion 1-2}), and  using the three term relation for
$\{\mathbb{Q}_n \}_{n \ge 0}$ and again (\ref{relacion 1-2}), we get
\begin{eqnarray*}
{A}_{n,i} \mathbb{P}_{n+1} &=& x_i \mathbb{Q}_{n}- \tilde{B}_{n,i}\mathbb{Q}_{n}-\tilde{C}_{n,i}\mathbb{Q}_{n-1}-A_{n,i} M_{n+1}\mathbb{P}_{n}\\
&=& x_i \mathbb{P}_{n}+ M_{n}x_i\mathbb{P}_{n-1}- \big(\tilde{B}_{n,i}+A_{n,i} M_{n+1}\big)\mathbb{P}_{n}\\
&~ &- \big( \tilde{C}_{n,i}+\tilde{B}_{n,i} M_{n}\big) \mathbb{P}_{n-1}- \tilde{C}_{n,i} M_{n-1} \mathbb{P}_{n-2},
\end{eqnarray*}
and by the induction hypothesis for  $x_i\mathbb{P}_{n-1}$, we obtain
\begin{eqnarray*}
A_{n,i} \mathbb{P}_{n+1} &=& x_i \mathbb{P}_{n} - \big(\tilde{B}_{n,i}-M_n A_{n-1,i}+A_{n,i} M_{n+1}\big)\mathbb{P}_{n}\\
& ~ & -\big(\tilde{C}_{n,i}-M_n {B}_{n-1,i}+\tilde{B}_{n,i} M_{n}\big)\mathbb{P}_{n-1}\\
& ~ & -\big(\tilde{C}_{n,i}M_{n-1}-M_{n}C_{n-1,i}\big)\mathbb{P}_{n-2}.
\end{eqnarray*}
Then taking into account (\ref{formula caracterizacion2}) we achieve the three term relation for $\{\mathbb{P}_n \}_{n \ge 0}$.

Also, we have
\begin{equation*}
\rank\, A_{n,i} = \rank\,\tilde{A}_{n,i} = \rank\,L_{n,i} = r_n^d, \quad 1\le i \le d,
\end{equation*}
and, for the joint matrix $A_n$, we  get
\begin{equation*}
\rank\, A_n = \rank\,\tilde{A}_n = \rank \, L_{n} = r_{n+1}^d .
\end{equation*}
To conclude, consider the linear functional $u$ defined by
$$\langle u, 1\rangle =1, \quad \langle u, \mathbb{P}_n\rangle =0, \quad n\ge 1,$$
which is well--defined since $\{\mathbb{P}_n : n\ge0\}$ is a basis of $\Pi^d$. We have just proved that $\{\mathbb{P}_n\}_{n\ge0}$ satisfies a three term relation (\ref{RR3TP}) and $A_n$ has full rank, then using the same arguments as in \cite[p. 80]{DX01}, we obtain that
\begin{equation}\label{aux}
\langle u, \mathbb{P}_k\,\mathbb{P}_j^t\rangle =0, \quad k\neq j.
\end{equation}

Next, we show that, for every $n \ge 0$, the symmetric and square matrix $H_n = \langle u, \mathbb{P}_n\,\mathbb{P}_n^t\rangle$ is invertible, that is, it has full rank. Taking into account \eqref{relacion 1-2} and \eqref{aux} we have
$$\langle u, \mathbb{Q}_n^t\rangle = 0, \quad n\ge2,$$
and expanding the linear functional $u$ in terms of the dual basis of  $\{\mathbb{Q}_n \}_{n \ge 0}$, and handling as in the proof of Lemma \ref{lema rango}, we can deduce that formula \eqref{M-H} holds, that is
$$M_n\, H_{n-1} = \tilde{H}_n\,(\sum_{i=1}^d a_i\,L^t_{n-1,i}), \quad n\ge 1.$$
We know that $\rank (\sum_{i=1}^d a_i\,L^t_{n-1,i}) = r^d_{n-1}$. Since the matrix $\tilde{H}_n$ is non--singular and the rank condition is invariant through non--singular matrices, we get
\begin{equation*}
\rank  \tilde{H}_n\,(\sum_{i=1}^d a_i\,L^t_{n-1,i}) = \rank (\sum_{i=1}^d a_i\,L^t_{n-1,i}) = r^d_{n-1},
\end{equation*}
and then
\begin{equation*}
 r^d_{n-1} = \rank M_n H_{n-1} \le \min\{\rank M_n, \rank H_{n-1}\} \le \rank H_{n-1} \le r_{n-1}^d.
\end{equation*}
Therefore $\rank H_{n-1} = r_{n-1}^d, \, n\ge 2$. Moreover, for $n=0$,
$$H_0 = \langle u, \mathbb{P}_0\,\mathbb{P}_0^t\rangle = \langle u, 1\rangle =1,$$
is an invertible matrix, and so for every $n \ge 0$, $H_n$ is invertible.
 Thus, $\{\mathbb{P}_n \}_{n \ge 0}$ is an OPS with respect to $u$ and the proof is completed.
\qquad $\square$

\bigskip

Now, we study the case when the monic polynomial system $\{\mathbb{P}_n \}_{n \ge 0}$ is {orthogonal}.

\begin{theorem} \label{caracterizacion ortogonalidad Q}
Let $\{\mathbb{P}_n \}_{n \ge 0}$ be
a monic orthogonal polynomial system satisfying the three term relation \eqref{RR3TP}.
Define the polynomial
system $\{\mathbb{Q}_n \}_{n \ge 0}$ by means of (\ref{relacion 1-2})
with  $\rank M_1=1$.

\noindent
Then $\{\mathbb{Q}_n\}_{n \ge 0}$ is a {monic orthogonal polynomial system} satisfying \eqref{RR3TQ} if and only if
formula (\ref{formula caracterizacion2}) holds and
\begin{eqnarray*}
 \rank\,  \tilde{C}_{n+1,i}&= & r_n^d,\quad  1 \le i \le d, \\
\rank\,  \tilde{C}_{n+1}^t &= & r_{n+1}^d,
\end{eqnarray*}
where
\begin{eqnarray*}
\tilde{A}_{n,i} &=& {A}_{n,i} = L_{n,i},\\
\tilde{B}_{n,i} &=& {B}_{n,i} + M_n A_{n-1,i}- \tilde{A}_{n,i} M_{n+1} , \\
\tilde{C}_{n,i} &=& C_{n,i} + M_n B_{n-1,i}- \tilde{B}_{n,i} M_{n} .
\end{eqnarray*}
\end{theorem}

\noindent
\textbf{Proof.} The necessary condition has already been proved in Theorem \ref{caracterizacion ortogonalidad 2}.

Conversely, writing (\ref{relacion 1-2}) for $n+1$, multiplying by $A_{n,i}$, and using the three term relation for $\mathbb{P}_n$, we get
\begin{equation*}
A_{n,i} \mathbb{Q}_{n+1} = x_i \mathbb{P}_{n}- \big( B_{n,i}-A_{n,i} M_{n+1}\big)\mathbb{P}_{n}-
C_{n,i} \mathbb{P}_{n-1}.
\end{equation*}
Using again (\ref{relacion 1-2}), we have
\begin{eqnarray*}
A_{n,i} \mathbb{Q}_{n+1} &=& x_i \mathbb{Q}_{n}- \big(B_{n,i}-A_{n,i} M_{n+1}\big)\mathbb{Q}_{n}\\
&~ &- M_{n} x_i \mathbb{P}_{n-1}- \big[ C_{n,i}-\big( B_{n,i}-A_{n,i} M_{n+1}\big) M_{n} \big] \mathbb{P}_{n-1}.
\end{eqnarray*}
Now, inserting (\ref{relacion 1-2}) in (\ref{RR3TP}), we obtain
\begin{equation*}
x_i \mathbb{P}_{n-1}=A_{n-1,i} \big(\mathbb{Q}_{n}- M_{n} \mathbb{P}_{n-1}\big)+ B_{n-1,i} \mathbb{P}_{n-1}+C_{n-1,i} \mathbb{P}_{n-2},
\end{equation*}
and therefore
\begin{eqnarray*}
\lefteqn{A_{n,i} \mathbb{Q}_{n+1}=x_i \mathbb{Q}_{n}- \big( B_{n,i}-A_{n,i} M_{n+1}+ M_{n}A_{n-1,i}\big)\mathbb{Q}_{n}}\\
&~ &-\Big[C_{n,i}+M_{n}B_{n-1,i}-M_{n}A_{n-1,i} M_{n}-\big(  B_{n,i}-A_{n,i} M_{n+1}\big)M_{n}\Big]\mathbb{P}_{n-1}\\
&~&-M_{n}C_{n-1,i}\mathbb{P}_{n-2}.
\end{eqnarray*}
In order to finish the proof, it is enough to replace $\mathbb{P}_{n-1}$ by $\mathbb{Q}_{n-1}- M_{n-1}\mathbb{P}_{n-2}$ and take into account the hypothesis (\ref{formula caracterizacion2}) and the expressions of $\tilde{A}_{n,i}$, $\tilde{B}_{n,i}$, and $\tilde{C}_{n,i}$. \qquad $\square$

\bigskip

It is worth to observe an essential difference between Theorems \ref{caracterizacion ortogonalidad 2} and  \ref {caracterizacion ortogonalidad Q}. In Theorem \ref{caracterizacion ortogonalidad 2}, starting from the orthogonality of $ \mathbb{Q}_{n}$, the conditions of full rank for the matrices $C_{n+1,i},\,  i=1, \dots , d$ and the joint matrix ${C}_{n+1}^t$  are deduced from Eq. (\ref{formula caracterizacion2}). However the situation is quite different if we assume the orthogonality of $ \mathbb{P}_{n}$. So in Theorem  \ref {caracterizacion ortogonalidad Q}, although the condition  which appears in the characterization of the orthogonality of $ \mathbb{Q}_{n}$ is the same Eq. (\ref{formula caracterizacion2}), it can not be deduced from it the requirements about the full rank of the matrices $\tilde{C}_{n+1,i},\, i=1, \dots, d$ and the joint matrix $\tilde{C}_{n+1}^t$.

Next, we give an example with $d=2$ to show that the required conditions of full rank for the corresponding matrices in Theorem \ref{caracterizacion ortogonalidad Q} are not superfluous. Indeed, for $ i=1,2$, consider the matrices
$$A_{n,i}= L_{n,i}, \quad C_{n,i}= -(L_{n-1,i})^t,  \quad  B_{n,i}=L_{n,i} C_{n+1,i} - C_{n,i} L_{n-1,i}. $$
{Observe that $B_{n,i}$ are $(n+1)\times (n+1)$ symmetric matrices with entries equal to  $0$ up to the entry $(n+1, n+1)$ of $B_{n,1}$, and the entry $(1, 1)$ of $B_{n,2}$ which are equal to $-1$. Obviously, $A_{n,i}$ and $C_{n+1,i}$ for $i=1, 2$, and  the joint matrices $A_n$ and $C_{n+1}^t$ have full rank. Then by Theorem \ref{Favard}, there exists a unique MOPS $\{\mathbb{P}_n\}_{n\geq0}$ with three term relation coefficients} $A_{n,i}, \, B_{n,i}$ and $C_{n,i}$, $i=1, 2.$

Consider $M_n = C_{n,1}$, and define a monic polynomial system $\{\mathbb{Q}_n \}_{n \ge 0}$ by
$$ \mathbb{Q}_n=\mathbb{P}_n + M_n \,\mathbb{P}_{n-1}, \quad n\ge 1.$$
Taking
\begin{eqnarray*}
\tilde{A}_{n,i} &=& {A}_{n,i}, \quad n\ge 0,\\
\tilde{B}_{n,i} &=& {B}_{n,i} + M_n A_{n-1,i}- \tilde{A}_{n,i} M_{n+1} , \quad n\ge 0,\\
\tilde{C}_{n,i} &=& C_{n,i} + M_n B_{n-1,i}- \tilde{B}_{n,i} M_{n} , \quad n\ge 1,
\end{eqnarray*}
straightforward computations lead to
\begin{eqnarray*}
\tilde{A}_{n,i} &=& {L}_{n,i}, \quad n\ge 0,\\
\tilde{B}_{n,1} &=& 0_{n+1,n+1} , \qquad  \tilde{B}_{n,2} ={B}_{n,2} , \quad n\ge 0,\\
\tilde{C}_{n,1} &=& {C}_{n,1}(I_n +  B_{n-1,1}) , \qquad \tilde{C}_{n,2} = {C}_{n,2}  , \quad n\ge 1.
\end{eqnarray*}
Moreover,
\begin{equation*}
\tilde{C}_{n,1} \, M_{n-1}={C}_{n,1} \,(I_n +  B_{n-1,1})\, C_{n-1,1}={C}_{n,1} \, C_{n-1,1}=M_n \, C_{n-1,1}, \quad n\ge 2,
\end{equation*}
and
\begin{equation*}
\tilde{C}_{n,2} \, M_{n-1}= M_n \, C_{n-1,2}, \quad n\ge 2.
\end{equation*}
Then by the previous results, the system $\{\mathbb{Q}_n\}_{n \ge 0}$ satisfies a three term relation with matrix coefficients $\tilde{A}_{n,i}, \, \tilde{B}_{n,i}$  and $\tilde{C}_{n,i}$, $i=1, 2.$ However $\rank\, \tilde{C}_{n,1}= n-1$, that is the matrix $\tilde{C}_{n,1}$ does  not have  full rank and therefore the system $\{\mathbb{Q}_n \}_{n \ge 0}$ is not orthogonal.

Note that concerning to the orthogonality of the linearly related polynomials, the above observation shows an important difference between the cases of several variables and one variable, and so the case of several variables is not a simple generalization of the case of one variable (see for instance \cite[Theorems 1 and 2]{MP95}).

\medskip

\section{Examples}

In this section we present several particular cases of orthogonal polynomial systems $\{\mathbb{P}_n \}_{n \ge 0}$ and $\{\mathbb{Q}_n \}_{n \ge 0}$ related by \eqref{relacion 1-2}, or equivalently \eqref{rel1-2non-monic}, giving the explicit expression of the involved matrices.

\subsection{Bivariate orthogonal polynomials related to Gaussian cubature formulas}

Linear combinations of orthogonal polynomials (quasi--orthogonal polynomials) in several variables have been considered  in connection with Gaussian cubature formulas. We apply our previous results to some examples developed  by Schmid and Xu (\cite{SX94}) based on some bivariate orthonormal polynomials introduced by Koornwinder in \cite{Koor75}:  Let $w(x)$ be a positive weight on an interval of $\mathbb{R}$. Let $\{p_n\}_{n \ge 0}$ be  the  sequence of orthonormal polynomials with respect to $w(x)$. It is well known that these polynomials  satisfies the three-term  recurrence formula
\begin{equation*}
x p_n(x)=a_n p_{n+1}(x)+ b_n p_n(x)+ a_{n-1} p_{n-1}(x), \quad n\ge 0\,,
\end{equation*}
where $p_0=1$ and $p_{-1}=0$.
Denote by $\{\mathbb{P}_n\}_{n \ge 0}$ the  sequence of bivariate orthonormal polynomials with respect to the weight function
$(u^2-4v)^{-1/2}W(u,v)$ where $W(u,v) = w(x) w(y)$, and  $u= x+y, \, v=xy$.

In \cite{SX94} the authors give an specific linear combination of the form $\mathbb{Q}_n  = \mathbb{P}_n  + M_{n,\rho} \, \mathbb{P}_{n-1}$ with
$$M_{n,\rho} = M_{n} =\rho \begin{pmatrix}
1 & \cdots & 0 & 0\\
\vdots & \ddots & \vdots & \vdots\\
0 &  \cdots & 1 & 0\\
0 &  \cdots & 0 & \sqrt{2}\\
0 &  \cdots & 0 & - \rho\\
\end{pmatrix}_{(n+1)\times n}, \quad \rho \in \mathbb R \backslash \{0\}, $$
in order to get  explicit Gaussian cubature formulas of degree $2n-2$.

Concerning to the orthogonality of the system $\{\mathbb{Q}_n\}_{n \ge 0}$, taking into account that $\{\mathbb{P}_n\}_{n \ge 0}$ is an  orthonormal polynomial system, Theorem \ref{caracterizacion ortogonalidad Q} yields the following characterization:

$\{\mathbb{Q}_n\}_{n \ge 0}$ is an orthogonal polynomial system if and only if
\begin{equation}\label{condicion Pn-ortonormales}
\tilde{C}_{n,i} \, M_{n-1} = M_n \, A_{n-2,i}^t, \quad n\ge 2,\quad i = 1, 2,
\end{equation}
and
\begin{eqnarray*}
\rank{\tilde{C}_{n,i}} &=& n, \quad n\ge 1,\quad i = 1, 2,\\
\rank{\tilde{C}_{n}^t} &=& n+1, \quad n\ge 1,\\
\end{eqnarray*}
with
\begin{eqnarray*}
\tilde{A}_{n,i}&=& {A}_{n,i}, \quad n\ge 0,\quad i = 1, 2,\\
\tilde{B}_{n,i} &=& {B}_{n,i} + M_n A_{n-1,i}- A_{n,i} M_{n+1}, \quad n\ge 1,\quad i = 1, 2,\\
\tilde{C}_{n,i} &=& A_{n-1,i}^t + M_n B_{n-1,i}- \tilde{B}_{n,i} M_n, \quad n\ge 1,\quad i = 1, 2.\\
\end{eqnarray*}

Using the explicit expressions for $M_n$  and for the matrices involved in the three--term relations satisfied by  $\mathbb{P}_n$  given in  \cite{SX94},  it is not too difficult to check that
$$\tilde{C}_{n,1} = \begin{pmatrix}
\lambda_{n,\rho} & 0 & \cdots & 0 & 0\\
0 & \lambda_{n,\rho} &  \cdots & 0 & 0\\
\vdots & \vdots & \ddots & \vdots & \vdots\\
0 & 0 & \cdots & \lambda_{n,\rho} & 0\\
0 & 0 & \cdots & \rho \,a_{n-2} & \sqrt{2}\,\lambda_{n,\rho}\\
0 & 0 & \cdots & -\sqrt{2}\,{\rho}^2\,a_{n-2} & -2\,\rho \,\lambda_{n,\rho}\\
\end{pmatrix}_{(n+1)\times n},\quad n \ge 2 \,,$$
with $$\lambda_{n,\rho}=a_{n-1}-{\rho}^2\,(a_{n-1}-a_n)+\rho\,(b_{n-1}-b_n).$$
Moreover, Eq. (\ref{condicion Pn-ortonormales}) for $i=1$ holds if and only if
 \begin{equation*}
 (a_{n-1}- a_{n-2})+ \rho (b_{n-1}- b_{n})+\rho^2 (a_{n}- a_{n-1})=0, \quad n\ge 2.
 \end{equation*}

 In particular, we analyze the orthogonality of the system $\{\mathbb{Q}_n\}_{n \ge 0}$  when $w(x)=w_{(\alpha, \beta)}(x)$ is a Chebyshev weight. Recall,

 \medskip

 a) Chebyshev of the first kind: $\alpha=\beta=-1/2$,\quad $w_{(-1/2,-1/2)}(x)=\frac{1}{\pi \sqrt{1-x^2}}$,
 $$a_0=1/\sqrt 2, \quad  a_n=1/2,\,\, n\ge1,  \quad {\rm{and}}\quad  b_n=0,\,\, n \ge 0.$$

 b) Chebyshev of the second kind: $\alpha=\beta=1/2$,\quad  $w_{(1/2,1/2)}(x)=\frac{2}{\pi} \,\sqrt{1-x^2}$,
 $$a_n=1/2,\,\, n\ge0, \quad {\rm{and}} \quad b_n=0,\,\, n \ge 0.$$

 c) Chebyshev of the third kind: $\alpha=-\beta=1/2$,\quad $w_{(1/2,-1/2)}(x)=\frac{1}{\pi}\,\sqrt{\frac {1-x}{1+x}}$,
 $$a_n=1/2,\,\, n\ge0,\quad b_0=-1/2  \quad {\rm{and}}\quad b_n=0,\,\, n \ge 1.$$

 d) Chebyshev of the fourth kind: $\alpha=-\beta=-1/2$,\quad $w_{(-1/2,1/2)}(x)=\frac{1}{\pi}\,\sqrt{\frac {1+x}{1-x}}$,
 $$a_n=1/2,\,\, n\ge0,\quad b_0=1/2  \quad {\rm{and}}\quad b_n=0,\,\, n \ge 1.$$

Observe that Eq. (\ref{condicion Pn-ortonormales}) for $i=1$ does not work  if $w(x)$ is the Chebyshev weight of the first kind
while for the remainder Chebyshev weights it holds. Thus, $\{\mathbb{Q}_n\}_{n \ge 0}$ is not an orthogonal polynomial system for $w(x)=w_{(-1/2,-1/2)}(x)$.

Moreover for the Chebyshev weights of the second, third and fourth kind, it is not difficult  to verify that Eq. (\ref{condicion Pn-ortonormales}) holds for $i=2$  since
the matrix $\tilde{C}_{n,2}$ takes the following form
$$\tilde{C}_{n,2} = \begin{pmatrix}
b_0\,a & a^2 & \cdots & 0 & 0 & 0\\
a^2 & 0 &  \cdots & 0 & 0 & 0\\
\vdots & \ddots & \ddots & \ddots & \vdots & \vdots\\
0 & 0 & \cdots & 0 & a^2 & 0\\
0 & 0 & \cdots & a^2 & \rho \, a^2 & \sqrt 2 \, a^2 \\
0 & 0 & \cdots & 0 & (1-\rho^2) \, a^2 & -\sqrt 2\, \rho \, a^2 \\
0 & 0 & \cdots & 0 & \sqrt 2\, \rho^3 \, a^2 & (1+2\rho^2) \, a^2\\
\end{pmatrix}_{(n+1)\times n},\quad n \ge 2, $$
where $a=1/2$.

Also, it is easy to check that  for $n\ge 2$  $\rank{\tilde{C}_{n,i}} = n, \,\, i = 1, 2$, and $\rank{\tilde{C}_{n}^t} = n+1$.

Finally, taking into account that for $n=1$ the expressions of the matrices $\tilde{C}_{1,i}, \, i=1,2$ are
$$\tilde{C}_{1,1} = \begin{pmatrix}
\sqrt 2 \,a + \sqrt 2 \,b_0 \, \rho\\
-2b_0\, \rho^2-2 \,a\,\rho\\
\end{pmatrix}_{2\times 1},$$
and
$$\tilde{C}_{1,2} = \begin{pmatrix}
\sqrt 2 \, a \, b_0 +\sqrt 2\, \rho(b_0^2-a^2)-\sqrt 2\,a \,b_0 \, \rho^2\\
a^2+ \rho^2(a-b_0^2)+b_0\, \rho^3\\
\end{pmatrix}_{2\times 1},$$
we have
$$\rank{\tilde{C}_{1,i}} = 1, \quad i = 1, 2, \quad \rm{and} \quad \rank{\tilde{C}_1^t} = 2,$$ for the following values of $\rho$ ($\rho \not=0$): For all values of $\rho$ in the case of Chebyshev of the second kind, for any $\rho\not=1$ in the case of third kind and for any $\rho\not=-1$ in the case of fourth kind.

Summarising, if $w(x)$ is a Chebyshev weight, the system $\{\mathbb{Q}_n\}_{n \ge 0}$ defined by
 \begin{equation*}
 \mathbb{Q}_n  = \mathbb{P}_n  + M_{n,\rho} \, \mathbb{P}_{n-1}
\end{equation*}
is an orthogonal polynomial system if and only if:

\medskip
a) $\rho = 0$ for Chebyshev  of the first kind,

 \medskip
b) $\rho \in \mathbb R $ for Chebyshev  of the second kind,

 \medskip
 c) $\rho  \in \mathbb R \backslash \{1\}$ for Chebyshev  of the third kind,

 \medskip
 d) $\rho \in \mathbb R  \backslash \{-1\}$ for Chebyshev   of the fourth kind.

\subsection{Koornwinder orthogonal polynomials}

We present some special examples of bivariate orthogonal polynomials generated by orthogonal polynomials of one variable satisfying a linear relation. To do this, we use the well known Koornwinder's method (\cite{DX01}, \cite{Koor75}).
More precisely, let $w_i(x)$, $i=1,2$, be two weight functions in one variable defined on the intervals $[a_i, b_i]$, respectively, and let $\rho(x)$ be a positive function in $[a_1, b_1]$ verifying either $\rho(x)$ is a polynomial of degree $1$ or $w_2(x)$ is a symmetric weight function and $\rho^2(x)$ is a polynomial of degree $\le 2$.

For $k\ge 0$, we denote by $\{q^{(k)}_{n}(x)\}_{n\ge 0}$  the sequence of univariate monic orthogonal polynomials with respect to the weight
function $\rho(x)^{2k+1}\,w_1(x)$, and by $\{r_n(y)\}_{n\ge 0}$ the sequence of monic orthogonal polynomials with respect to
$w_2(y)$. Consider the polynomials of two
variables of total degree $n$ given by
$$
Q_{n-k,k}(x,y)=q^{(k)}_{n-k}(x)\, \rho(x)^k\, r_k
\Big(\frac{y}{\rho(x)}\Big), \qquad 0\le k\le n,
$$
which are orthogonal  with respect to the weight function
$$W(x,y) =w_1(x)w_2\Big(\frac{y}{\rho(x)}\Big),$$
on the region $\{(x,y)\in\mathbb{R}^2: a_1 < x < b_1, a_2\,\rho(x) < y < b_2\,\rho(x)\}$ (see \cite{DX01}, p. 55).

Suppose that there exists a sequence of monic polynomials in one variable $\{p_{n}(x)\}_{n\ge 0}$ orthogonal with respect to the weight $\widetilde{w}_1(x)$ satisfying the relation
$$(x-\xi) \, w_1(x) = \widetilde{w}_1(x),$$
where $\xi\in \mathbb{R}\setminus [a_1, b_1]$ is a fixed real number. Then {for $k\ge 0$,} the two following weight functions
\begin{eqnarray*}
w_1^{(k)}(x) &=&  \rho(x)^{2k+1}\,w_1(x) ,\\
\widetilde{w}_1^{(k)}(x) &=& \rho(x)^{2k+1}\,\widetilde{w}_1(x),
\end{eqnarray*}
are related by
$$(x-\xi) \, w_1^{(k)}(x) = \widetilde{w}_1^{(k)}(x).$$

Therefore (see \cite{MP95}), for a fixed $k\ge0$, the polynomials $\{q^{(k)}_{n}(x)\}_{n\ge 0}$ and the polynomials $\{p^{(k)}_{n}(x)\}_{n\ge 0}$ orthogonal with respect to the weight function $\widetilde{w}_1^{(k)}(x)$, satisfy the linear relation
\begin{equation*}
q^{(k)}_n(x) = p^{(k)}_n(x) + a^{(k)}_n\, p^{(k)}_{n-1}(x), \quad n\ge 1,
\end{equation*}
where $a^{(k)}_n\in \mathbb{R}$, for $ k\ge 0$.

Using this fact, the orthogonal polynomials in two variables $Q_{n-k,k}(x,y)$ satisfy the linear relation
$$Q_{n-k,k}(x,y) = P_{n-k,k}(x,y) + a^{(k)}_{n-k} \,P_{n-1-k,k}(x,y), $$
where
$$P_{n-k,k}(x,y) = p^{(k)}_{n-k}(x)\, \rho(x)^k\, r_k
\Big(\frac{y}{\rho(x)}\Big)$$
are bivariate polynomials orthogonal with respect to the weight function
$$\widetilde{W}(x,y) = \lambda(x,y)\, W(x,y),$$
with $\lambda(x,y)=(x-\xi)$. In this way, the orthogonal polynomial systems
\begin{eqnarray*}
&~&\{\mathbb{P}_n(x,y)\}_{n\ge0} = \{(P_{n,0}(x,y), P_{n-1,1}(x,y),\ldots, P_{0,n}(x,y))^t\}_{n\ge 0},\\
&~&\{\mathbb{Q}_n(x,y)\}_{n\ge0} = \{(Q_{n,0}(x,y), Q_{n-1,1}(x,y),\ldots, Q_{0,n}(x,y))^t\}_{n\ge 0},
\end{eqnarray*}
satisfy the matrix linear relation
\begin{equation}\label{relacion-matricial}
\mathbb{Q}_n(x,y) = \mathbb{P}_n (x,y)+ M_n \, \mathbb{P}_{n-1}(x,y),
\end{equation}
where $M_n$ is given by
$$M_n = \begin{pmatrix}
a_n^{(0)} & 0 & \cdots & 0\\
0 & a_{n-1}^{(1)} &  \cdots & 0\\
\vdots & \vdots & \ddots & \vdots\\
0 & 0 & \cdots & a_1^{(n-1)}\\
0 & 0 & \cdots & 0
\end{pmatrix}_{(n+1)\times n}.$$

\bigskip

Using this procedure, we can deduce relations between some families of
well known orthogonal polynomials in two variables. As far as we know, these relations are new.

\subsubsection{Orthogonal polynomials on the unit disk}

Orthogonal polynomials in two variables on the unit disk $B^2=\{(x,y)\in \mathbb{R}^2: x^2 + y^2\le 1\},$ (the so--called {\it disk polynomials})
are associated with the inner product
$$(f,g)_\mu = c_\mu\int_{B^2} f(x,y)g(x,y) W^{(\mu)}(x,y) dx\,dy,$$
where
$$W^{(\mu)}(x,y) = (1-x^2-y^2)^\mu, \qquad \mu > -1,$$
is the weight function, and $c_\mu$ is the normalization constant in order to have $(1,1)_\mu =1$.

Using Koornwinder's tools, disk polynomials can be defined from Jacobi polynomials as
$$Q^{(\mu)}_{n-k,k}(x,y) = P^{(\mu+\frac{1}{2}+k,\mu+\frac{1}{2}+k)}_{n-k}(x)\,(1-x^2)^{\frac{k}{2}}\,P^{(\mu,\mu)}_{k}((1-x^2)^{-\frac{1}{2}}\,y), \quad 0\le k\le n,$$
taking
\begin{eqnarray*}
w_1(x) &=& (1-x^2)^{\mu},\quad x\in [-1,1], \quad \mu >-1,\\
 w_2(y) &=& (1-y^2)^\mu,\quad y\in [-1,1], \quad \mu >-1,\\
\rho(x) &=& \sqrt{1-x^2}.
\end{eqnarray*}
Since monic Jabobi polynomials satisfy the relation (see \cite[Chapter 22]{AS72})
$$P^{(\alpha,\alpha)}_{n}(x)= P^{(\alpha+1,\alpha)}_{n}(x) - \frac{n}{2n+2\alpha+1} \,P^{(\alpha+1,\alpha)}_{n-1}(x),$$
we can write
$$
 Q^{(\mu)}_{n-k,k}(x,y) = P^{(\mu+1)}_{n-k,k}(x,y) - \frac{n-k}{2n+2\mu+2} \,P^{(\mu+1)}_{n-1-k,k}(x,y),
$$
where
$$P^{(\mu+1)}_{n-k,k}(x,y) = P^{(\mu+\frac{3}{2}+k, \mu+\frac{1}{2}+k)}_{n-k}(x) \,(1-x^2)^{\frac{k}{2}}\,P^{(\mu,\mu)}_{k}((1-x^2)^{-\frac{1}{2}}\,y)
$$
are Koornwinder polynomials associated with
the weight function on the unit disk
$$\widetilde{W}^{(\mu)}(x,y) = (1-x)\,W^{(\mu)}(x,y).$$
Then relation (\ref{relacion-matricial}) holds, where the matrix $M_n$ is given by
$$M_n = \frac{-1}{2n+2\mu+2}\,\begin{pmatrix}
n & 0 & \cdots & 0\\
0 & n-1 &  \cdots & 0\\
\vdots & \vdots & \ddots & \vdots\\
0 & 0 & \cdots & 1\\
0 & 0 & \cdots & 0
\end{pmatrix}_{(n+1)\times n}.$$

\subsection{Tensor product of polynomials in one variable}

When $\rho(x)=1$, Koornwinder's method leads to tensor product of orthogonal polynomials in one variable. This case can be rewritten for moment functionals in the following way. Let  $v_x$ and $w_y$ be two quasi--definite moment functionals  (acting on variables $x$ and $y$, respectively) and let $\{q_n(x)\}_{n\ge 0}$ and $\{r_n(y)\}_{n\ge 0}$ be their respective sequences of  orthogonal polynomials in the variables $x$ and $y$.

The polynomials
$$Q_{n-k,k}(x,y) = q_{n-k}(x) \, r_{k}(y), \quad 0\le k \le n,$$
are orthogonal with respect to the {\it composition} moment functional $v$
$$
\langle v, f(x,y) \rangle :=\langle v_{x}, \langle w_{y}, f(x, y)\rangle\rangle = \langle w_{y}, \langle v_{x}, f(x, y)\rangle\rangle, \quad \forall f\in \Pi^2,
$$
namely $v=v_{x}\circ w_{y}= w_{y}\circ v_{x}$.

Suppose that there exists a quasi--definite moment functional $u_x$ related with $v_x$ by $$(x-\xi)\,v_x=u_x,$$ where $\xi$ is a fixed real number. Then the orthogonal polynomials  $\{q_{n}(x)\}_{n\ge 0}$ are linearly related with the monic orthogonal polynomials
$\{p_{n}(x)\}_{n\ge 0}$ associated with the quasi--definite moment functional $u_x$. In this way (\cite{MP95}), there exist non zero constants $\{a_n\}_{n\ge 1}$ such that
\begin{equation*}
q_n(x) = p_n(x) + a_n\, p_{n-1}(x), \quad n\ge 1.
\end{equation*}
Then the polynomials $Q_{n-k,k}(x,y)$ satisfy a linear relation
of the form
$$Q_{n-k,k}(x,y) = P_{n-k,k}(x,y) + a_{n-k} \,P_{n-1-k,k}(x,y), \quad n\ge 1, \quad 0\le k\le n-1,$$
where
$$P_{n-k,k}(x,y) = p_{n-k}(x)\, r_k(y),$$
are bivariate orthogonal polynomials associated with the moment functional
$u = u_x \circ w_{y}$.

Moreover, both moment functionals are related by
$$\lambda(x,y)\,v=u,$$
with $\lambda(x,y)=(x-\xi)$, and besides relation (\ref{relacion-matricial}) holds where the matrices $M_n$ are given by
\begin{equation}\label{matriz tensor}
M_n = \begin{pmatrix}
a_n & 0 & \cdots & 0\\
0 & a_{n-1} &  \cdots & 0\\
\vdots & \vdots & \ddots & \vdots\\
0 & 0 & \cdots & a_1\\
0 & 0 & \cdots & 0
\end{pmatrix}_{(n+1)\times n}.
\end{equation}
Obviously, an analogous situation occurs when $\{r_k(y)\}_{k\ge 0}$ is {\it linearly related} with other orthogonal polynomial sequence.

Next we present two new examples of orthogonal polynomial systems in two variables generated from tensor products of polynomials in one variable which are orthogonal with respect to quasi--definite moment  functionals.

\subsubsection{Krall Laguerre--Laguerre orthogonal polynomials}

Consider the classical Laguerre moment functional
$$
u_x = x^{\alpha}\,e^{-x}, \quad \alpha >-1,
$$
and the modification given by
$$
v_x = x^{-1} \, u_x +\frac{\Gamma(\alpha+1)}{\alpha+1-a_1}\,\delta_0,
$$
where $a_1 \not=0$ is a real constant such that $\alpha+1-a_1\not=0$.

Recall that the action of the functional $(x-c)^{-1} \, u$ over a polynomial is defined by (see \cite{Mar91})
$$\langle (x-c)^{-1} \, u, p \rangle := \langle u, \frac{p(x)-p(c)}{x-c} \rangle.$$

The moment functional $v_x$ is quasi--definite if and only if (see \cite[p. 896]{AMPR11}) either
\begin{eqnarray*}
\alpha_n:=\Gamma(n)\,\Gamma(\alpha+1)\,(\alpha+1-a_1)+(a_1-1)\,\Gamma(n+\alpha)&\not= &0, \quad n\ge 2, \quad \text{for} \quad \alpha \not=0,
\end{eqnarray*} or
\begin{eqnarray*}
\widetilde{\alpha}_n:=(a_1-1)\,(1+1/2+\dots +1/(n-1))+1 &\not=& 0, \quad n\ge 2, \quad \text{for} \quad \alpha =0.
\end{eqnarray*}

Let $\{L_n^{(\alpha)}\}_{n\ge0}$ and $\{Q_n\}_{n\ge0}$ be the sequences of monic polynomials orthogonal with respect to the functionals $u_x$ and $v_x$, respectively. Since the following relation
$$xv_x=u_x$$
holds, we have
$$Q_n(x)=L_n^{(\alpha)}(x)+a_n\,L_{n-1}^{(\alpha)}(x), \quad n\ge1.$$
In  \cite{AMPR11}, it was obtained the explicit expression of the coefficients $a_n, \, n\geq2,$
\begin{eqnarray*}
 a_n=\begin{cases} \displaystyle \quad \frac{\alpha_{n+1}}{\alpha_n}, \quad \alpha \not=0,\\
 \\
 \displaystyle n\,\frac{\widetilde{\alpha}_{n+1}}{\widetilde{\alpha}_n},  \quad \alpha =0. \end{cases}
\end{eqnarray*}

 Let $w_y$ be any quasi--definite moment functional, and define
$$v = v_x\circ w_y, \quad u=u_x\circ w_y.$$
These moment functionals satisfy the relation $\lambda(x,y)v=u$ where $\lambda(x,y)=x$. If both moment funcionals are quasi--definite, the respective bivariate orthogonal polynomials satisfy the relation (\ref{relacion-matricial}) and the matrix $M_n$ is given by (\ref{matriz tensor}).

\subsubsection{Krall Jacobi--Jacobi orthogonal polynomials}

In Section 4 of \cite{APPR12}, the authors consider the classical Jacobi moment functional
$$u_x = (1-x)^{\alpha}(1+x)^{\beta}, \quad \alpha, \beta >-1,$$
and the modification
$$
v_x = (1-x)^{-1} \, u_x +\langle u_x,1 \rangle\,\frac{\alpha+\beta+2}{2(\alpha+1)+a_1(\alpha+\beta+2)}\, \delta_1,
$$
where $\langle u_x,1 \rangle=\int_{-1}^1 (1-x)^{\alpha}(1+x)^{\beta}\,dx$, and $a_1\not=0$ is a parameter satisfying
$$2(\alpha+1)+a_1(\alpha+\beta+2) \not=0.$$

As it was proven in \cite{APPR12}, the moment functional $v_x$ is quasi--definite if and only if either
$$
\alpha_n:=\Gamma(\alpha+1)\Gamma(\alpha+\beta+2)\Gamma(n)\Gamma(n+\beta)
+M \Gamma(\beta+1)\Gamma(n+\alpha)\Gamma(n+\alpha+\beta)\not=0,\quad  n\ge2 ,
$$
for $\alpha \not=0$, and $
M:=-\frac{\displaystyle{2(\beta+1)+a_1(\alpha+\beta+1)(\alpha+\beta+2)}}{\displaystyle{2(\alpha+1)+a_1(\alpha+\beta+2)}},$
or
$$\widetilde{\alpha}_n:=\displaystyle \frac{2 (\beta+2)}{2+a_1(\beta+2)}-(\beta+1) \sum_{i=1}^{n-1}\left(\frac{1}{i}+\frac{1}{\beta+i}\right)\not=0,\quad  n\ge2, \quad \text{for} \quad \alpha=0.
$$
Let $\{P_n^{(\alpha, \beta)}\}_{n\ge0}$ and $\{Q_n\}_{n\ge0}$ be the sequences of monic polynomials orthogonal with respect to the functionals $u_x$ and $v_x$, respectively. Since
$$(1-x)\,v_x = u_x$$
then
$$Q_n(x)=P_n^{(\alpha, \beta)}(x)+a_n \,P_{n-1}^{(\alpha, \beta)}(x), \quad n\ge1. $$
In \cite{APPR12}, the authors give an explicit expression of the parameters $a_n, \, n \geq 2,$ in terms of the free parameter $a_1$
\begin{eqnarray*}
a_n=\begin{cases} \displaystyle \frac{-2}{(2n+\alpha+\beta)(2n+\alpha+\beta-1)} \frac{\alpha_{n+1}}{\alpha_n},\quad &  \alpha \not=0,\\
\\
\quad \displaystyle  \frac{-2n(n+\beta)}{(2n+\beta)(2n+\beta-1)} \frac{\widetilde{\alpha}_{n+1}}{\widetilde{\alpha}_n},\quad  & \alpha =0. \end{cases}
\end{eqnarray*}

Then for any quasi--definite moment functional $w_y$,  the two quasi--definite moment functionals $v=v_x\circ w_y$ and $u=u_x\circ w_y$  satisfy the relation
$$\lambda(x,y)v=u,$$
where $\lambda(x,y)=1-x$. Thus, the bivariate orthogonal polynomials associated with $u$ and $v$ satisfy the relation (\ref{relacion-matricial}), where the matrix $M_n$ is given explicitly by (\ref{matriz tensor}).

\subsection{Adjacent families of classical orthogonal polynomials in several variables}

This subsection is devoted to deduce {\it relations between adjacent families} of classical orthogonal polynomials in several variables, that is, to give some polynomials as linear combinations of polynomials of the same family with different values of their parameters. These relations can be seen as a generalization of the ones for Jacobi and Laguerre polynomials in one variable.

\subsubsection{Classical orthogonal polynomials on the simplex (Appell polynomials)}

Classical polynomials on the simplex (see \cite{DX01}, p. 46) are orthogonal with respect to the
inner product
$$(f,g)_\kappa = \omega_\kappa\int_{T^d} f(\mathtt{x}) g(\mathtt{x}) W^{(\kappa)}(\mathtt{x}) d\mathtt{x},$$
on the simplex in $\mathbb{R}^d$,
$$
T^d = \{\mathtt{x} = (x_1, x_2, \ldots,x_d) \in \mathbb{R}^d:\,x_1, x_2, \ldots, x_d \geq 0,\quad 1-|\mathtt{x}|_1\geq 0\},
$$
where the weight function is given by
$$W^{(\kappa)}(\mathtt{x}) = x_1^{\kappa_1-1/2} \,x_2^{\kappa_2-1/2}\cdots x_d^{\kappa_d-1/2}\,(1-|\mathtt{x}|_1)^{\kappa_{d+1}-1/2},\quad \kappa_i > -\frac{1}{2},$$
for $\mathtt{x}\in T^d$ and $|\mathtt{x}|_1 = x_1 + \cdots + x_d$ is the usual $\ell^1$ norm.
Denoting $\kappa = (\kappa_1,\kappa_2, \ldots, \kappa_{d+1})$, the normalization constant $\omega_\kappa$ is taken in order to have $(1,1)_\kappa=1$, and it is given by
$$\omega_\kappa= \frac{\Gamma(|\kappa|+\frac{d+1}{2})}{\Gamma(\kappa_1+\frac{1}{2})\,\cdots\,\Gamma(\kappa_{d+1} +\frac{1}{2})},$$
with $|\kappa|= \kappa_1 + \cdots + \kappa_{d+1}$.

\noindent
We will use the following notation. For $\mathtt{x} = (x_1, \ldots, x_d) \in \mathbb{R}^d$, we define the {\it truncation} of $\mathtt{x}$ as
$$\mathtt{x}_0 =0, \quad \mathtt{x}_j = (x_1, \ldots, x_j), \quad 1\le j \le d.$$
Observe that $\mathtt{x}_d = \mathtt{x}$. Associated with $\nu=(\nu_1,\dots,\nu_d)$, we define
$$\nu^j=(\nu_j,\dots,\nu_d), \quad 1\le j \le d.$$
Moreover, we denote by $e_1 = (1,0, \ldots,0)$ the first vector of the canonical basis.

For a multi--index $\nu=(\nu_1,\dots,\nu_d)\in \mathbb{N}^d_0$, a basis of orthonormal polynomials on the simplex is given by (\cite[p. 47]{DX01}\footnote{The formula that appears in this Subsection has been rewritten using the document published by the author in \tt{http://pages.uoregon.edu/yuan/paper/Errata.pdf} })
\begin{eqnarray}\label{base-orth-triang}
P_{\nu}^{(\kappa)}(\mathtt{x}) &=& [h_\nu^{(\kappa)}]^{-1} \, \prod_{j=1}^d\left(\frac{1-|\mathtt{x}_j|_1}{1-|\mathtt{x}_{j-1}|_1}\right)^{|\nu^{j+1}|}\,p^{(a_j, b_j)}_{\nu_j}\left(\frac{2\,x_j}{1-|\mathtt{x}_{j-1}|_1} -1\right)\nonumber\\
&=& [h_\nu^{(\kappa)}]^{-1} (1-x_1)^{|\nu^2|}p^{(a_1, b_1)}_{\nu_1}(2\,x_1 -1)\nonumber\\
&~& \times \prod_{j=2}^d\left(\frac{1-|\mathtt{x}_j|_1}{1-|\mathtt{x}_{j-1}|_1}\right)^{|\nu^{j+1}|}
\,p^{(a_j, b_j)}_{\nu_j}\left(\frac{2\,x_j}{1-|\mathtt{x}_{j-1}|_1} -1\right),
\end{eqnarray}
where $a_j = |\kappa^{j+1}|+2|\nu^{j+1}| + \frac{d-j-1}{2}$, $b_j = \kappa_j - \frac{1}{2}$, the polynomials $\{p^{(a, b)}_m(t)\}_{m\ge 0}$ are the orthonormal Jacobi polynomials in one variable,
$$
[h_{\nu}^{(\kappa)}]^{2} = \frac{\prod_{j=1}^d (|\kappa^j| + 2|\nu^{j+1}| + \frac{d-j+2}{2})_{2\nu_j}}{(|\kappa| + \frac{d+1}{2})_{2|\nu|}},
$$
and $(a)_n = a (a+1)\cdots (a+n-1)$ denotes the usual Pochhammer symbol for $a\in \mathbb{R}$ and $n\ge 0$, with the convention $(a)_0=1$.

The following relation between orthonormal families of Jacobi polynomials with different parameters can be easily deduced from formula $(22.7.19)$ in \cite{AS72}:

\begin{equation}\label{lin_rel_jac_orth}
p_m^{(a,b)}(t) = c_m^{(a,b)}\,p_m^{(a,b+1)}(t) +  d_m^{(a,b)}\,p_{m-1}^{(a,b+1)}(t),\quad m\ge0
\end{equation}
where
$$
c_m^{(a,b)}= \left[\frac{2(m+b+1)(m+a+b+1)}{(2m+a+b+2)(2m+a+b+1)}\right]^{1/2}, \,\,
d_m^{(a,b)}= \left[\frac{2\,m(m+a)}{(2m+a+b+1)(2m+a+b)}\right]^{1/2}.
$$
Then, substituting in \eqref{base-orth-triang}, we get
\begin{eqnarray}
P_{\nu}^{(\kappa)}(\mathtt{x}) &=& [h_{\nu}^{(\kappa)}]^{-1} (1-x_1)^{|\nu^2|}\left[c_{\nu_1}^{(a_1,b_1)}p^{(a_1,b_1+1)}_{\nu_1}(2x_1-1) +
d_{\nu_1}^{(a_1,b_1)}p^{(a_1, b_1+1)}_{\nu_1-1}(2x_1-1)\right]\nonumber\\
&~&\times \prod_{j=2}^d\left(\frac{1-|\mathtt{x}_j|_1}{1-|\mathtt{x}_{j-1}|_1}\right)^{|\nu^{j+1}|}
\,p^{(a_j,b_j)}_{\nu_j}\left(\frac{2\,x_j}{1-|\mathtt{x}_{j-1}|_1} -1\right)\nonumber\\
&=& \frac{h_{\nu}^{(\kappa+e_1)}}{h_{\nu}^{(\kappa)}}\, c_{\nu_1}^{(a_1,b_1)}\,P_{\nu}^{(\kappa+e_1)}(\mathtt{x}) +
 \frac{h_{\nu-e_1}^{(\kappa+e_1)}}{h_{\nu}^{(\kappa)}}\, d_{\nu_1}^{(a_1,b_1)}P_{\nu-e_1}^{(\kappa+e_1)}(\mathtt{x}),\label{1-2}
\end{eqnarray}
where the second summand vanishes for $\nu_1=0$. Observe that $\{P_{\nu}^{(\kappa+e_1)}(\mathtt{x}): |\nu| =n, n\ge 0\}$ are the orthonormal polynomials defined by \eqref{base-orth-triang}, associated
with the inner product
$$(f,g)_{\kappa+e_1} = \omega_{\kappa+e_1}\int_{T^d} f(\mathtt{x}) g(\mathtt{x}) W^{(\kappa+e_1)}(\mathtt{x}) d\mathtt{x},$$
and
$$W^{(\kappa+e_1)}(\mathtt{x}) =x_1\,W^{(\kappa)}(\mathtt{x})= x_1^{\kappa_1+1/2}x_2^{\kappa_2-1/2}\cdots x_d^{\kappa_d-1/2}\, (1-|\mathtt{x}|_1)^{\kappa_{d+1}-1/2}.$$
Next, we represent relation \eqref{1-2} in matrix form like \eqref{rel1-2non-monic} using the orthonormal polynomial systems on the simplex $\{\mathbb{P}_n^{(\kappa)}\}_{n\ge0}$ and $\{\mathbb{P}_n^{(\kappa+e_1)}\}_{n\ge0}$.

Let $n\ge1$, and let $\alpha_1, \alpha_2, \ldots, \alpha_{r_n^d}$ be the elements in
$\{\alpha \in \mathbb{N}^d_0: |\alpha| = n\}$ arranged according to the
reverse lexicographical order. We will denote the components of the multi--index
$\alpha_i$ as
$$\alpha_i = (\alpha_{i,1}, \alpha_{i,2}, \ldots, \alpha_{i,d}), \qquad i=1, 2, \ldots, r_n^d.$$
Observe that $\alpha_{i,1}\ge 1$ for $i=1, 2, \ldots, r_{n-1}^d $ and $\alpha_{i,1}=0$ for $r^d_{n-1}< i \le r^d_n$. Moreover
$\alpha_1-e_1, \alpha_2-e_1, \ldots, \alpha_{r^d_{n-1}}-e_1$ are the elements in
$\{\beta \in \mathbb{N}^d_0: |\beta| = n-1\}$ arranged again according to the reverse lexicographical order.

For $n\ge 1 $, we define the matrices
$$\hat{K}_n^{(1)}=\left(\begin{array}{cccccc}
{c}_{\alpha_1}^{(1)}&                    &      &                              &        & \bigcirc  \\
                    &{c}_{\alpha_2}^{(1)}&      &                              &        &         \\
                    &                    &\ddots&                              &        &         \\
                    &                    &      &{c}_{\alpha_{r^d_{n-1}}}^{(1)}&        &         \\
                    &                    &      &                              & \ddots &         \\
\bigcirc            &                    &      &                              & &{c}_{\alpha_{r_n^d}}^{(1)}
\end{array}\right), \quad
\hat{M}_n^{(1)}=\left(\begin{array}{cccc}
{d}_{\alpha_1}^{(1)}&                   &       & \bigcirc                  \\
                   &{d}_{\alpha_2}^{(1)}&       &                           \\
                   &                   &\ddots &                           \\
\bigcirc           &                   &       &{d}_{\alpha_{r^d_{n-1}}}^{(1)}\\
\hline
                   &                   &       &                             \\
\bigcirc           & \cdots            &\cdots &\bigcirc                      \\
                   &                   &       &                               \\
\end{array}\right),$$
of respective sizes $r_n^d\times r_n^d$ and $r_n^d\times r^d_{n-1}$, where
$$
{c}_{\alpha_i}^{(1)} = \frac{h_{\alpha_i}^{(\kappa+e_1)}}{h_{\alpha_i}^{(\kappa)}}\, c_{\alpha_{i,1}}^{(a_1,b_1)}, \quad i=1,2,\ldots,r^d_n, \qquad {d}_{\alpha_i}^{(1)} = \frac{h_{\alpha_i-e_1}^{(\kappa+e_1)}}{h_{\alpha_i}^{(\kappa)}}\, d_{\alpha_{i,1}}^{(a_1,b_1)}, \quad i=1, 2, \ldots, r^d_{n-1}.
$$
Then, \eqref{1-2} reads as
\begin{equation*}
\mathbb{P}_n^{(\kappa)}(\mathtt{x}) = \hat{K}_n^{(1)}\,\mathbb{P}_n^{(\kappa+e_1)}(\mathtt{x}) + \hat{M}_n^{(1)} \, \mathbb{P}_{n-1}^{(\kappa+e_1)}(\mathtt{x}), \quad n\ge 1.
\end{equation*}
Notice that $\hat{K}_n^{(1)}$ is non singular and $\hat{M}_n^{(1)}$ has full rank.

Likewise we could have replaced formula \eqref{lin_rel_jac_orth} in \eqref{base-orth-triang} for every Jacobi polynomial $p^{(a_j, b_j)}_{\nu_j}(t)$, for $1\le j\le d$ fixed. Then, a similar procedure shows that
\begin{equation*}
\mathbb{P}_n^{(\kappa)}(\mathtt{x}) = \hat{K}_n^{(j)}\,\mathbb{P}_n^{(\kappa+e_j)}(\mathtt{x}) + \hat{M}_n^{(j)} \, \mathbb{P}_{n-1}^{(\kappa+e_j)}(\mathtt{x}), \quad n\ge 1,
\end{equation*}
holds for the matrices
$$\hat{K}_n^{(j)}= \hbox{\rm{diag}}\{{c}_{\alpha_i}^{(j)}: i=1, 2, \ldots, r^d_n\}, \qquad \hat{M}_n^{(j)}= L_{n-1,j}^t\, \hbox{\rm{diag}}\{{d}_{\alpha_i}^{(j)}: i=1, 2, \ldots, r^d_n,\,\text{s.t.}\, \alpha_{i,j}\ge 1 \},$$
with
$$
{c}_{\alpha_i}^{(j)} = \frac{h_{\alpha_i}^{(\kappa+e_j)}}{h_{\alpha_i}^{(\kappa)}}\, c_{\alpha_{i,j}}^{(a_j,b_j)}, \qquad {d}_{\alpha_i}^{(j)} = \frac{h_{\alpha_i-e_j}^{(\kappa+e_j)}}{h_{\alpha_i}^{(\kappa)}}\, d_{\alpha_{i,j}}^{(a_j,b_j)},
$$
where $e_j = (0, \ldots, 0, 1, 0, \ldots, 0)$ is the jth vector of the canonical basis.

\subsubsection{Multiple Jacobi polynomials on the $d$--cube}

Multiple Jacobi polynomials on the cube (\cite[p. 37]{DX01}) are orthogonal with respect to the multiple Jacobi weight function
$$W^{(a,b)}_J(\mathtt{x}) = \prod_{i=1}^d \, (1-x_i)^{a_i}\, (1+x_i)^{b_i},$$
on the cube $[-1,1]^d$ of $\mathbb{R}^d$, where $\mathtt{x}=(x_1, \ldots, x_d)$, and
$$a=(a_1, \ldots, a_d), \qquad b=(b_1,\ldots, b_d), \qquad a_i, b_i >-1.$$
According to \cite{DX01}, an orthogonal basis is given in terms of standard Jacobi polynomials by
\begin{equation}\label{mult-Jac}
P_\nu(\mathtt{x}; W^{(a,b)}_J) = P_{\nu_1}^{(a_1,b_1)}(x_1) \cdots P_{\nu_d}^{(a_d,b_d)}(x_d), \quad |\nu|=n.
\end{equation}
The following relations between adjacent families of Jacobi polynomials can be found in
(\cite[Chapter 22]{AS72}):
\begin{equation}\label{ra}
P_m^{(a,b)}(t) = f_m^{(a,b)} \, P_m^{(a+1,b)}(t) - g_m^{(a,b)}\,P_{m-1}^{(a+1,b)}(t),\quad m\ge0
\end{equation}
\begin{equation}\label{rb}
P_m^{(a,b)}(t) = f_m^{(a,b)} \, P_m^{(a,b+1)}(t) + g_m^{(b,a)}\,P_{m-1}^{(a,b+1)}(t),\quad m\ge0
\end{equation}
where
$$
f_m^{(a,b)} = \frac{m+a+b+1}{2\,m + a + b+1}, \quad g_m^{(a,b)} = \frac{m+b}{2\,m + a + b+1}.
$$
Let $j$ be fixed with $1\le j\le d$. Then, we can substitute \eqref{ra} in \eqref{mult-Jac}, and we obtain
\begin{eqnarray*}
P_\nu(\mathtt{x}; W^{(a,b)}_J) &=& f_{\nu_j}^{(a_j,b_j)} \,P_{\nu_1}^{(a_1,b_1)}(x_1) \cdots P_{\nu_j}^{(a_j+1,b_j)}(x_j)\cdots P_{\nu_d}^{(a_d,b_d)}(x_d) \\
&~& - g_{\nu_j}^{(a_j,b_j)}\, P_{\nu_1}^{(a_1,b_1)}(x_1) \cdots P_{\nu_j-1}^{(a_j+1,b_j)}(x_j)\cdots P_{\nu_d}^{(a_d,b_d)}(x_d) \\
&=& f_{\nu_j}^{(a_j,b_j)}\,P_\nu(\mathtt{x}; W^{(a+e_j,b)}_J) - g_{\nu_j}^{(a_j,b_j)}\,P_{\nu-e_j}(\mathtt{x}; W^{(a+e_j,b)}_J).
\end{eqnarray*}

As in the above example, we denote by $\alpha_1, \alpha_2, \ldots, \alpha_{r_n^d}$  the elements in
$\{\alpha \in \mathbb{N}^d_0: |\alpha| = n\}$ arranged according to the reverse lexicographical order with $\alpha_i = (\alpha_{i,1}, \alpha_{i,2}, \ldots, \alpha_{i,d})$ for $i=1, 2, \ldots, r_n^d$, and by $e_j = (0, \ldots, 0, 1, 0, \ldots, 0)$ the jth vector of the canonical basis.

In this way, if $\{\mathbb{P}_n(\mathtt{x};W^{(a,b)}_J)\}_{n\ge 0}$ denotes the classical Jacobi polynomial system on the cube defined as above, for $1\le j \le d$, and $n\ge 1$, we get
\begin{equation*}
\mathbb{P}_n(\mathtt{x};W_J^{(a,b)}) = \hat{K}_n^{(j)}(a,b)\,\mathbb{P}_n(\mathtt{x};W^{(a+e_j,b)}_J) - \hat{M}_n^{(j)}(a,b) \, \mathbb{P}_{n-1}(\mathtt{x};W^{(a+e_j,b)}_J),
\end{equation*}
where
\begin{eqnarray*}
\hat{K}_n^{(j)}(a,b)&=& \hbox{\rm{diag}}\{{f}_{\alpha_{i,j}}^{(a_j,b_j)}: i=1, 2, \ldots, r^d_n\},\\ \hat{M}_n^{(j)}(a,b) &=& L_{n-1,j}^t\, \hbox{\rm{diag}}\{{g}_{\alpha_{i,j}}^{(a_j,b_j)}: i=1, 2, \ldots, r^d_n,\,\text{ s.t.} \, \alpha_{i,j} \ge1 \}.
\end{eqnarray*}

In a similar way, using relation (\ref{rb}) we obtain
\begin{equation*}
\mathbb{P}_n(\mathtt{x};W^{(a,b)}_J) = \hat{K}_n^{(j)}(a,b)\,\mathbb{P}_n(\mathtt{x};W^{(a,b+e_j)}_J) + \hat{M}_n^{(j)}(b,a) \, \mathbb{P}_{n-1}(\mathtt{x};W^{(a,b+e_j)}_J),
\end{equation*}
where
$$\hat{M}_n^{(j)}(b,a) = L_{n-1,j}^t\, \hbox{\rm{diag}}\{{g}_{\alpha_{i,j}}^{(b_j,a_j)}; i=1, 2, \ldots, r^d_n,\,\text{s.t.} \, \alpha_{i,j}\ge 1 \}.$$

\subsubsection{Multiple Laguerre polynomials on $\mathbb{R}^d_+$}

Multiple Laguerre polynomials are orthogonal with respect to the weight function (\cite[p. 51]{DX01})
$$W^{(\kappa)}_L(\mathtt{x}) = \mathtt{x}^\kappa\, e^{-|\mathtt{x}|_1}, \qquad \mathtt{x}\in \mathbb{R}^d_+,$$
which is the product of Laguerre weights in one variable.

As in the above case, multiple Laguerre polynomials defined by
$$P_\nu(\mathtt{x};W^{(\kappa)}_L) = L^{(\kappa_1)}_{\nu_1}(x_1)\cdots L^{(\kappa_d)}_{\nu_d}(x_d),\quad |\nu| =n, $$
form a mutually orthogonal basis associated with $W^{(\kappa)}_L$.

Using formula (5.1.13) in \cite{Sz78},
$$L^{(a)}_m(t) = L^{(a+1)}_m(t) - L^{(a+1)}_{m-1}(t),\quad m \ge 0$$
we get the following relation between adjacent families of Laguerre polynomials
$$P_\nu(\mathtt{x};W^{(\kappa)}_L) = P_\nu(\mathtt{x};W^{(\kappa+e_j)}_L) - P_{\nu-e_j}(\mathtt{x};W^{(\kappa+e_j)}_L),$$
and $e_j = (0, \ldots, 0, 1, 0, \ldots, 0)$ is the jth vector of the canonical basis.
In a matricial form, we express above relation as
$$\mathbb{P}_n(\mathtt{x};W^{(\kappa)}_L) = \mathbb{P}_n(\mathtt{x};W^{(\kappa+e_j)}_L) - L^t_{n-1,j} \, \mathbb{P}_{n-1}(\mathtt{x};W^{(\kappa+e_j)}_L),$$
where $\{\mathbb{P}_n(\mathtt{x};W^{(\kappa)}_L)\}_{n\ge 0}$ denotes the Laguerre polynomial system on $\mathbb{R}^d_+$.

\bigskip

\textbf{Acknowledgements} The authors would like to thank   the referee for his valuable suggestions and  comments  which helped to improve the paper.

\end{document}